\documentclass{amsart}

\usepackage{amssymb}

\usepackage[all]{xy}

\newtheorem{thm}{Theorem}%[section]

\newtheorem{lem}[thm]{Lemma}
\newtheorem{cor}[thm]{Corollary}

\newtheorem{prop}[thm]{Proposition}

\newtheorem{conj}[thm]{Conjecture}
\theoremstyle{definition}
\newtheorem{defn}[thm]{Definition}

\newtheorem{say}[thm]{}
\newtheorem{exmp}[thm]{Example}

\newtheorem{ques}[thm]{Question}    %!!!!!!!!!!!!!!!!!!!!

\newtheorem*{ack}{Acknowledgments}      % \renewcommand{\theack}{} 

\newtheorem{defn-thm}[thm]{Definition--Theorem}  %!!!!!!!!!!!!!!!!!!!!!!!!
\newtheorem{defn-lem}[thm]{Definition--Lemma}  %!!!!!!!!!!!!!!!!!!!!!!!!
\theoremstyle{remark}

\renewcommand{\c}[0]{{\mathbb C}}  

\renewcommand{\o}[0]{{\mathcal O}} 
\newcommand{\z}[0]{{\mathbb Z}}

\renewcommand{\a}[0]{{\mathbb A}}

\newcommand{\p}[0]{{\mathbb P}}

\newcommand{\q}[0]{{\mathbb Q}}
\newcommand{\map}[0]{\dasharrow}
\newcommand{\qtq}[1]{\quad\mbox{#1}\quad}

\newcommand{\pic}[0]{\operatorname{Pic}}
\newcommand{\gal}[0]{\operatorname{Gal}}

\newcommand{\rank}[0]{\operatorname{rank}}

\newcommand{\codim}[0]{\operatorname{codim}}

\newcommand{\Hom}[0]{\operatorname{Hom}}

\newcommand{\aut}[0]{\operatorname{Aut}}

\newcommand{\sing}[0]{\operatorname{Sing}}

\newcommand{\chow}[0]{\operatorname{Chow}}

\newcommand{\cl}[0]{\operatorname{Cl}}

\newcommand{\onto}[0]{\twoheadrightarrow}

\newcommand{\simq}[0]{\sim_{\q}}

\newcommand{\tsum}[0]{\textstyle{\sum}}

\newcommand{\shom}[0]{\operatorname{\mathcal{H}\!\it{om}}}

\newcommand{\simb}[0]{\stackrel{bir}{\sim}}

\def\into{\DOTSB\lhook\joinrel\to}

\def\loccoh#1.#2.#3.#4.{H^{#1}_{#2}(#3,#4)}

\DeclareMathAlphabet{\mathchanc}{OT1}{pzc}%
                                {m}{it}

\newcommand{\alb}[0]{\operatorname{Alb}}
\newcommand{\albr}[0]{\operatorname{Alb}^{\rm rat}}
\begin{document}
\bibliographystyle{amsalpha}

%\today

\title{Deformations of elliptic Calabi--Yau manifolds}
\author{J\'anos Koll\'ar}

\maketitle

The aim of this note is to answer some questions 
about Calabi--Yau manifolds that were raised during the
workshop {\it String Theory for Mathematicians,}
which was held  at the Simons Center for Geometry and Physics.

F-theory posits that the ``hidden dimensions'' constitute 
a Calabi--Yau 4-fold $X$ that has an elliptic structure with a section.
That is, there are morphisms $g:X\to B$ whose general fibers are
elliptic curves and $\sigma:B\to X$ such that
$g\circ \sigma=1_B$; see \cite{Vafa96,1998hep.th.2093D}.  In his lecture 
Donagi asked the following.

\begin{ques}\label{ques1}
 Is every small deformation of an elliptic Calabi--Yau manifold
 also  an elliptic Calabi--Yau manifold?
\end{ques}

\begin{ques} \label{ques2}
Is there a good numerical characterization of
elliptic Calabi--Yau manifolds?
\end{ques}

Note that a good answer to Question \ref{ques2} should
give an answer to Question \ref{ques1}.
The answers to these problems  are quite sensitive to 
which variant of the definition of  Calabi--Yau manifolds one uses.
For instance,  a general deformation of 
 $(\mbox{Abelian variety})\times(\mbox{elliptic curve})$
 has no elliptic fiber space structure and
every elliptic K3 surface has non-elliptic deformations.
We prove  in Section 5 that these are essentially the only
such examples, even for singular Calabi--Yau varieties
(\ref{thm.ques1.general}). In the smooth case, the answer is
especially simple.

\begin{thm}\label{thm.ques1} Let $X$ be an elliptic Calabi--Yau
manifold such that $H^2(X, \o_X)=0$. Then 
every small  deformation of $X$ is also an 
elliptic Calabi--Yau manifold.
\end{thm}

In dimension 3 this was proved in \cite{MR1314743, MR1490200}.

Our results on Question \ref{ques2} are less complete.
Let  $L_B\in H^2(B,\q)$ be an ample cohomology class and set
$L:=g^*L_B$. 
We interpret Question \ref{ques2} to mean: {\it Characterize
pairs $(X, L)$ that are elliptic fiber spaces.}
Following \cite{Wilson-MR1010159, oguiso93}, one is led to the following.

\begin{conj} \label{ques4} A Calabi--Yau manifold $X$ is elliptic iff there is a
$(1,1)$-class $L\in H^2(X, \q)$ such that 
$(L\cdot C)\geq 0$ for every algebraic curve $C\subset X$,
$\bigl(L^{\dim X}\bigr)=0$
and $\bigl(L^{\dim X-1}\bigr)\neq 0$.
% as an element of $H^{2\dim X-2}\bigl(X, \q\bigr)$.
\end{conj}

%Conjecture \ref{ques4} should hold if $K_X\simq 0$.
%and $H^1(X, \o_X)=0$.
For threefolds, the more general results of \cite{oguiso93, MR1314743}
imply  Conjecture \ref{ques4}  if $L$ is effective or
$\bigl(L\cdot c_2(X)\bigr)\neq 0$.

As in \cite{Wilson-MR1010159, oguiso93, MR1314743}, in higher dimensions 
we  study the interrelation of
$L$ and of  the
second  Chern class $c_2(X)$. 
By a result of \cite{MR949837}  
$\bigl(L^{n-2}\cdot c_2(X)\bigr)\geq 0$ 
and we distinguish two (overlapping) cases.
\begin{itemize}
\item (Main case) If $\bigl(L^{n-2}\cdot c_2(X)\bigr)> 0$ then 
 Conjecture \ref{ques4} is solved in
(\ref{main.char.thm.cor.q2}). 
We also check that all  elliptic   Calabi--Yau manifolds
with a section belong to this class (\ref{isotriv.section.few}).
\item (Isotrivial case) These are the examples where
$X\to B$ is an analytically
 locally trivial fiber bundle over a dense open subset
of $B$. An explicit construction, up-to birational equivalence,
is given in (\ref{isortiv.standard.form.exmp})
but I do not have a numerical characterization.
\end{itemize}

Following \cite{oguiso93} and \cite[Lect.10]{MP97}, 
the plan is to put both questions in the more general framework
of the Abundance Conjecture \cite[4.6]{Reid83};
see (\ref{ques3}--\ref{charact.conj}) for the precise formulation.

This approach suggests that the key  is to understand
the rate of growth of $h^0(X, L^m)$.
If  $(X,L)$ is elliptic, then $h^0(X, L^m)$ grows like
$m^{\dim X-1}$.
% and the hardest part is to show that
% $h^0(X, L^m)$ has this expected order of growth.
Given a pair $(X,L)$, the most important 
deformation-invariant quantity is the holomorphic
Euler characteristic  
$$
\chi(X, L^m)=h^0(X, L^m)-h^1(X, L^m)+h^2(X, L^m)\cdots
$$
The difficulty is that in our case 
$h^0(X, L^m)$ and $h^1(X, L^m) $ both grow like $m^{\dim X-1}$
and they cancel each other out. That is
$$
\chi(X, L^m)=O\bigl(m^{\dim X-2}\bigr).
$$
For the main series, $\chi(X, L^m) $ does grow like
$m^{\dim X-2} $ which implies that $h^0(X, L^m)$ grows at least like
$m^{\dim X-2}$.

For the isotrivial series the  order of growth
of  $\chi(X, L^m) $ is even smaller; in fact $\chi(X, L^m) $ can 
be identically zero.
However, if $(X,L)$ is elliptic, this happens only if a finite cover of $X$ is
birational to a product, so these are not particularly interesting
examples.

Several  of the ideas of this paper can be traced back to other sources.
 Sections 2--4 owe a lot to \cite{Kawamata85b, oguiso93, MR1314743, MR2779478};
Sections 5--6 to \cite{MR0417458, k-lars};
 Sections 7--8 to \cite{k-shaf, nakayama}
and to some  old results of Matsusaka. 
Ultimately the origin of many of these methods is in  the work
of Kodaira on elliptic surfaces \cite[Sec.12]{MR0184257}.
(See  \cite[Secs.V.7--13]{bpv} for a more modern treatment.)

\section{Calabi--Yau fiber spaces}

For many reasons it is of interest to study proper morphisms
with connected fibers $g':X'\to B$
whose general fibers are birational to
Calabi--Yau varieties. A special case of the Minimal Model Conjecture,
proved by \cite{lai-2009, hacon-xu},  implies that
every such fiber space is birational to a
projective morphism
with connected fibers $g:X\to B$ where $X$ has terminal
singularities and its canonical class $K_X$ is
relatively trivial, at least rationally. That is,
there is a Cartier divisor $F$ on $B$ such that
$mK_X\sim g^*F$ for some $m>0$.

We will work with varieties with log terminal singularities,
or  later even with klt pairs $(X, \Delta)$ but
I will state the main results for smooth varieties as well.
See \cite[Sec.2.3]{km-book} for the definitions and basic properties
of the singularities we use.
Note also that, even if one is primarily interested in  smooth  Calabi--Yau 
 varieties $X$, the natural setting is to allow
at least  canonical  singularities on $X$
and at least log terminal singularities on the base $B$ of the
elliptic fibration.

\begin{defn}\label{Calabi--Yau.fib.defn}
In this paper a  {\it  Calabi--Yau variety} is a
projective variety $X$ with log terminal singularities
such that $K_{X}\simq 0$, that is, $mK_X$ is linearly equivalent to
0 for some $m>0$. Note that by \cite{kaw} this is equivalent to
assuming that $(K_X\cdot C)=0$ for every curve $C\subset X$.

Note that we allow a rather broad definition of
Calabi--Yau varieties. This is very  natural for 
algebraic geometry but less so for physical considerations.

 A {\it  Calabi--Yau fiber space} is a
proper morphisms
with connected fibers  $g:X\to B$
onto a normal variety
where $X $ has log terminal (or possibly log canonical) singularities and
 $K_{X_g}\simq 0$ where $X_g\subset X$ is a general fiber.

We say that $g:X\to B$  is an {\it elliptic}
(or {\it Abelian} or ...) fiber space if in addition
 general fibers are elliptic curves (or Abelian varieties or ...).
Our main  interest is in the elliptic case, but  in Sections 7--8
we also study  the general setting.

Let $X$ be a projective,  log terminal variety
 and $L$ a $\q$-Cartier $\q$-divisor
(or divisor class) on $X$. We say that  $(X, L)$
is a  Calabi--Yau fiber space if there is a
 Calabi--Yau fiber space $g:X\to B$
and an ample $\q$-Cartier $\q$-divisor $L_B$ on $B$ 
such that $L\simq g^*L_B$. 

In general, a divisor $L$ is called {\it semi-ample} 
if it is the pull-back of an ample divisor by a morphism and 
{\it nef} if $(L\cdot C)\geq 0$ for every irreducible curve $C\subset X$.
Every semi-ample  divisor is nef but the converse usually fails.
However, the hope is that for Calabi--Yau varieties nef and  semi-ample are 
equivalent;
see (\ref{ques3}--\ref{charact.conj}).

We say that a  Calabi--Yau fiber space $g:X\to B$ is
{\it relatively minimal} if 
$K_X\simq g^*F$ for some
$\q$-Cartier $\q$-divisor $F$ on $B$.
This condition is automatic if $X$ itself is Calabi--Yau.
(These are called crepant log structures in \cite{kk-singbook}.)

If $g:X\to B$ is a relatively minimal  Calabi--Yau fiber space
and $X$ has canonical (resp.\ log terminal) singularities
then every other  relatively minimal  Calabi--Yau fiber space
$g':X'\to B$ that is birational to $X$ also has
canonical (resp.\ log terminal) singularities.

By \cite{naka88}, if $X$ has  log terminal singularities
then $B$ has rational  singularities, more precisely,
there is an effective divisor $D_B$ such that $(B, D_B)$ is
klt.
\end{defn}

\begin{say}[Elliptic  threefolds]\label{ell.Calabi--Yau.3f.say}
Elliptic  threefolds have been studied in detail. The papers
 \cite{Wilson-MR1010159, MR1109635, MR1260943, MR1217381, MR1242006, MR1272699, MR1272978, MR1314743, MR1401771, MR1918053, MR1929795,  2010arXiv1008.2018C, hul-klo, 2011arXiv1107.2043K} 
give rather complete descriptions of
 their local and global structure.
However, neither of Questions \ref{ques1}--\ref{ques2} was
fully answered for threefolds.

By contrast, not even the  local structure of elliptic  fourfolds
is understood.  Double covers of the $\p^1$-contractions
described in 
\cite{MR1620110} give some rather surprising examples; there are
probably much more complicated ones as well.
\end{say}

\begin{defn}\label{excpet.div.defn}
Let $g:X\to B$ be a  morphism between normal varieties.
A divisor $D\subset X$ is called 
{\it horizontal} if $g(D)=B$, 
{\it vertical} if $g(D)\subset B$ has codimension  $\geq 1$  and
{\it exceptional}
if $g(D)$ has codimension $\geq 2$
in $B$. 

If $g$ is birational, the latter coincides with the usual notion of
exceptional divisors but the above version makes sense even if
$\dim X>\dim B$.
(If $g$ is birational then there are no horizontal divisors, so this notion is
not used in that case.)
\end{defn}

\begin{say} \label{various.say}
We see in (\ref{nef.vertucal.div.lem}) that 
if $X$ is smooth (or $\q$-factorial), $g$ is a Calabi--Yau fiber space and 
 $D\subset X$ is exceptional then 
 $D$
is not $g$-nef. Thus, by \cite{lai-2009, hacon-xu} the 
$ (X,\epsilon D)$ Minimal Model Program  over $B$ 
(cf. \cite[Sec.3.7]{km-book}) contracts $D$. 
Thus every Calabi--Yau fiber space $g_2:X_2\to B_2$ is
birational to a relatively minimal 
 Calabi--Yau fiber space $g_1:X_1\to B_1=B_2$ that has no exceptional divisors.
Furthermore, again using \cite{lai-2009, hacon-xu}  and applyting 
(\ref{kaw.abound}) it is also birational to  a 
Calabi--Yau fiber space $g:X\to B$ where $B$ is also 
$\q$-factorial.
(In general  $g:X\to B$  is not unique.)
Thus, in birational geometry,  it is reasonable to focus on the study of 
 relatively minimal 
 Calabi--Yau fiber spaces  $g:X\to B$  without  exceptional divisors
where both $X$ and $B$ are $\q$-factorial and log terminal.

From the point of view of  $F$-theory  it is especially
interesting to study the examples
$g':X'\to B$ with a section $\sigma':B\to X'$
where $X'$ itself is Calabi--Yau. 
In this case the so called {\it Weierstrass model} is a 
%one obtains a distinguished
 relatively minimal model  without  exceptional divisors
that can be explicitly constructed 
as follows.

%One obtains the Weierstrass model as follows.
 Let $L_B$ be an ample divisor on $B$.
Then $\sigma'(B)+m{g'}^*L_B$ is nef and big on $X'$,
hence a large multiple of it is base point free
(cf.\ \cite[Sec.3.2]{km-book}). 
This gives a morphism $h:X'\to X$ where
$X$ is still Calabi--Yau (usually with canonical singularities)
and $g:X\to B$ has a section $\sigma:B\to X$
whose image is $g$-ample. Thus every fiber of $g$
has dimension $1$ and so $g:X\to B$  has no exceptional divisors.

Furthermore, $R^1h_*\o_{X'}=0$ which implies that every deformation of
$X'$ comes from a deformation of $X$; see (\ref{deform.conj.say}).
\end{say}

The next result says  that once  $g:X\to B$ looks like a
  relatively minimal Calabi--Yau fiber space
outside a subset of codimension $\geq 2$ then it is a
 relatively minimal Calabi--Yau fiber space.

\begin{prop} \label{codim2.conj}
Let  $g:X\to B$ be a projective fiber space
with $X $ log terminal.
Assume the following.
\begin{enumerate}
\item There are no $g$-exceptional divisors (\ref{excpet.div.defn}).
\item  There is a closed subset $Z\subset B$ of codimension
$\geq 2$ such that $K_X$ is numerically trivial on the fibers
over $B\setminus Z$.
\item $B$ is $\q$-factorial.
\end{enumerate}
Then $g:X\to B$ is a relatively minimal Calabi--Yau fiber space.
\end{prop}

Proof. First note that, as a very special case of  (\ref{kaw.abound}),
there is a  $\q$-Cartier $\q$-divisor $F_1$ on  $B\setminus Z$
such that 
$$
K_X|_{X\setminus g^{-1}(Z)}\simq g^* F_1.
$$
Since $B$ is $\q$-factorial, $F_1$ extends to a
$\q$-Cartier $\q$-divisor $F$ on  $B$.

Thus every point $b\in B$ has an open neighborhood
$b\in U_b\subset B$ and an integer $m_b>0$ such that
$$
\o_X\bigl(m_bK_X\bigr)|_{g^{-1}(U_b\setminus Z)}\cong g^* 
\o_{U_b}\bigl(m_bF|_{U_b}\bigr)\cong g^*  \o_{U_b}
\cong \o_{g^{-1}(U_b)}.
$$
By (1),  $g^{-1}(Z)$ has codimension $\geq 2$ in 
$g^{-1}(U_b)$ and hence the constant 1 section of
$\o_{g^{-1}(U_b\setminus Z)}$ extends to a global section
of $\o_X\bigl(m_bK_X\bigr)|_{g^{-1}(U_b)}$
that has neither poles not zeros.
Thus 
$$
\o_X\bigl(m_bK_X\bigr)|_{g^{-1}(U_b)}\cong \o_{g^{-1}(U_b)}.
$$
Since this holds for every $b\in B$, we conclude that
$K_X\simq g^*F$.\qed

\section{The main case}

The next theorem gives a characterization of
the main series of elliptic Calabi--Yau fiber spaces.
(For the log version see (\ref{logell.conj.say}).)
The proof is quite short but it relies on
auxiliary results that are proved in the next two sections.

\begin{thm}\label{main.char.thm.td2}
 Let $X$ be a projective variety of dimension $n$
with log terminal singularities
and $L$ a Cartier divisor on $X$. 
Assume that $K_X$ is nef and 
$\bigl(L^{n-2}\cdot \operatorname{td}_2(X)\bigr)>0$
where $\operatorname{td}_2(X)$ is the second Todd class of $X$
(\ref{RR.ratsing.say}).
Then  
$(X,L)$ is a relatively minimal,   elliptic fiber space
iff
\begin{enumerate}
\item $L$ is nef,
\item $L-\epsilon K_X$ is nef for $0\leq \epsilon \ll 1$,
\item $(L^n)=(L^{n-1}\cdot K_X)=0$ and
\item $(L^{n-1})$ is nonzero in $H^{2n-2}(X,\q)$.
\end{enumerate}
\end{thm}

Note that if 
$(X,L)$ is a relatively minimal,   elliptic fiber space
then $L$ is semi-ample (\ref{Calabi--Yau.fib.defn}) 
and, as we see in (\ref{pf.main.char.thm.td2}),
the only hard part of (\ref{main.char.thm.td2}) is to show that
conditions (\ref{main.char.thm.td2}.1--4) imply that $L$ is semi-ample.
In particular,  (\ref{main.char.thm.td2})
also holds over fields that are not algebraicaly closed.

This immediately yields the following partial
answer to Question \ref{ques1}. 

\begin{cor}\label{main.char.thm.cor.q2}
  Let $X$ be a smooth, projective  variety of dimension $n$
and $L$ a Cartier divisor on $X$. 
Assume that  $K_X\simq 0$ and 
$\bigl(L^{n-2}\cdot c_2(X)\bigr)>0$.
Then  
$(X,L)$ is  an   elliptic fiber space
iff
\begin{enumerate}
\item $L$ is nef,
\item $(L^n)=0$ and
\item $(L^{n-1})$ is nonzero in $H^{2n-2}(X,\q)$. 
\end{enumerate}
\end{cor}

\begin{defn} \label{kappa.nu.defn}
Let $Y$ be a projective variety and $D$ a Cartier divisor on $X$.
If $m>0$ is sufficiently divisible, then, up-to birational equivalence,
 the map given by global sections of  $\o_Y(mD)$
$$
 Y\map  I(Y,D)\simb I_m(Y,D)\into \p\bigl(H^0(Y, \o_Y(mD))\bigr) \qtq{is independent of $m$.}
$$
It is called the {\it Iitaka fibration} of $(Y,D)$.
The {\it Kodaira dimension} of $D$ (or of $(Y,D)$)
is $\kappa(D)=\kappa(Y,D):=\dim I(Y,D)$.

If $D$ is nef, the
{\it numerical dimension}  of $D$ (or of $(Y,D)$),
denoted by $\nu(D)$ or $\nu(Y,D)$, is the largest
natural number $r$ such that the self-intersection
$(D^r)\in H^{2r}(Y,\q)$ is nonzero. Equivalently,
$(D^r\cdot H^{n-r})>0$ for some (or every)
ample divisor $H$.

It is easy to see that $\kappa(D)\leq \nu(D)$. 
This was probably first observed by Matsusaka as
a corollary of his theory of variable intersection cycles; see
\cite{Matsusaka72} or \cite[p.515]{MR0379494}.

\end{defn}

\begin{say}[Proof of (\ref{main.char.thm.td2})] 
\label{pf.main.char.thm.td2}
First note that $\kappa(L)\geq n-2$ by 
(\ref{there.are.sections.lem}). 
We will also need this for some perturbations of $L$.

Set $L_m:=L-\frac1{m} K_X$.
For $m\gg 1$ we see that $L_m$ is nef,
$\bigl(L_m^{n-2}\cdot \operatorname{td}_2(X)\bigr)>0$
and $(L_m^{n-1})$ is nonzero in $H^{2n-2}(X,\q)$.
 Note that $mL=K_X+mL_m$ hence
$$
m^n(L^n)=\tsum_{i=0}^n\, m^{n-i}\bigl(K_X^i\cdot L_m^{n-i}\bigr).
%\qtq{and}
%m^{n-1}(L^{n-1}\cdot K_X)=
%\sum_{i=0}^{n-1} m^{n-1-i}\bigl(K_X^{i+1}\cdot L_m^{n-1-i}\bigr).
$$
Since $K_X$ and $L_m$ are both nef, all the terms on the right hand side are
$\geq 0$. Their sum is zero by assumption, hence
$\bigl(K_X^i\cdot L_m^{n-i}\bigr)=0$ for every $i$. 
 Thus (\ref{there.are.sections.lem})
also applies to $L_m$ and we get  that
$\kappa(L_m)\geq n-2$. 

We can now apply (\ref{oguiso.prop}) with
$\Delta=0$ and $D:=2mL_m$ 
and $K_X+2mL_m=2mL_{2m}$ 
to conclude that
$\nu(L_m)\leq \kappa(L_{2m})$.  Since we know that
$\nu(L_m)=\dim X-1$ we conclude that
$\kappa(L_{2m})=\dim X-1$. 

Finally use (\ref{kaw.abound})
with $S=(\mbox{point})$, $2mL$ instead of $L$ and $a=1$
to obtain that some multiple of $L$ is semi-ample.
That is, there is a morphism with connected fibers $g:X\to B$
and an ample $\q$-divisor $L_B$ such that $L\simq g^*L_B$. 
Note that  $\bigl(L^{\dim B}\bigr)\neq 0$ but
 $\bigl(L^{\dim B+1}\bigr)=0$ so comparing with (3--4) we see that
$\dim B=\dim X-1$. 
By the adjunction formula, the canonical class of the general fiber
is proportional to $(L^{n-1}\cdot K_X)=0$, thus
$g:X\to B$ is an elliptic fiber space.
\qed 
\end{say}

We have used the following theorem
due to \cite{Kawamata85b} and \cite{MR2779478}.
%see also  \cite[6-1-11]{ka-ma-ma}.

\begin{thm} \label{kaw.abound}
Let $(X,\Delta)$ be an irreducible, projective, klt pair and
$g:X\to S$ a morphism with generic fiber $X_g$.
Let $L$ be a $\q$-Cartier $\q$-divisor on $X$. Assume that
\begin{enumerate}
\item $L$ and $L-K_X-\Delta$ are $g$-nef and
\item  $\nu\bigl((L-K_X-\Delta)|_{X_g}\bigr)=
\kappa\bigl((L-K_X-\Delta)|_{X_g}\bigr)=
\nu\bigl(((1+a)L-K_X-\Delta)|_{X_g}\bigr)=
\kappa\bigl(((1+a)L-K_X-\Delta)|_{X_g}\bigr)$
for some $a>0$.
\end{enumerate}
Then there is a factorization
$g:X\stackrel{h}{\to} B\stackrel{\pi}{\to} S$ and a
$\pi$-ample  $\q$-Cartier $\q$-divisor $L_B$ on $B$
such that $L\simq h^*L_B$.\qed
\end{thm}

\section{Adjoint systems of large Kodaira dimension}

The following is modeled on \cite[2.4]{oguiso93}.

\begin{prop} \label{oguiso.prop}
 Let $(X,\Delta)$ be a projective, klt pair 
such that $K_X+\Delta$ is pseudo-effective,
that is, its cohomology class is a limit of effective classes.
Let  $D$ be
an effective,  nef, $\q$-Cartier $\q$-divisor on $X$ such
 that  $\kappa(K_X+\Delta+D)\geq \dim X-2$. 

Then  $\nu(D)\leq \kappa(K_X+\Delta+D)$.
\end{prop}

Proof. There is nothing to prove if 
$\kappa(K_X+\Delta+D)= \dim X$.
Thus assume that  $\kappa(K_X+\Delta+D)\leq  \dim X-1$ and
let $g:X\map B$ be the Iitaka fibration (cf.\ \cite[2.1.33]{laz-book}).
After some blow-ups we may assume that $g$ is a morphism and $X, B$ are smooth.

The generic fiber of $g$ is a smooth curve or surface
$(S,\Delta_S)$ such  that $K_S+\Delta_S$
 is pseudo-effective.
Since abundance holds for curves and surfaces \cite[Sec.11]{k-etal},
this implies that $\kappa(K_S+\Delta_S)\geq 0$.
Furthermore,  by Iitaka's theorem (cf.\ \cite[2.1.33]{laz-book})
$\kappa(K_S+\Delta_S+D|_S)= 0$.

If $D$ is disjoint from $S$ then, by (\ref{nef.vertucal.div.cor}.2),
$\nu(D)\leq \dim B=\kappa(K_X+\Delta+D)$ and we are done.
Otherwise $D|_S$ is an effective,  nonzero, nef divisor on $S$.
We obtain a contradiction by proving that
$\kappa(K_S+\Delta_S+D|_S)\geq 1$.

If $S$ is a curve, then $\deg D|_S>0$ hence
$\kappa(K_S+\Delta_S+D|_S)\geq \kappa(D|_S)=1$.
If $S$ is a surface, then $\kappa(K_S+\Delta_S+D|_S)\geq 1$
is proved in (\ref{surf.nu=1.lem}). \qed

\begin{lem} \label{surf.nu=1.lem} Let $(S,\Delta_S)$ be a projective,
klt surface such that $\kappa(K_S+\Delta_S)\geq 0$.
Let $D$ be a nonzero, effective, nef $\q$-divisor. Then
$\kappa(K_S+\Delta_S+D)\geq 1$.
\end{lem}

Proof. Since $\kappa(K_S+\Delta_S+D)\geq \kappa(K_S+\Delta_S)$
we only need to consider the case when $\kappa(K_S+\Delta_S)=0$.
Let $\pi:(S,\Delta_S)\to (S^m,\Delta_S^m)$ be the  minimal model.
It is obtained by repeatedly contracting curves that have
negative intersection number with $K_S+\Delta_S $.
These curves also have negative intersection number with 
$K_S+\Delta_S +\epsilon D$
for $0<\epsilon\ll 1$. Thus 
$$
\pi: (S,\Delta_S+\epsilon D)\to (S^m,\Delta_S^m+\epsilon D^m)
$$
is also  the minimal model 
and $(S^m,\Delta^m+\epsilon D^m) $ is klt
for $0<\epsilon\ll 1$.
By the Hodge index theorem,
every effective divisor contracted by $\pi$ has negative
self-intersection, thus $D^m$ is again 
 a nonzero, effective, nef $\q$-divisor. 

Since abundance holds for klt surface pairs 
(cf.\ \cite[Sec.11]{k-etal}),
we see that  
$K_{S^m}+\Delta^m\simq 0$ and 
$\kappa\bigl( K_{S^m}+\Delta^m+\epsilon D^m\bigr)\geq 1$.
Since $D$ is effective, we obtain that
$
\kappa(K_S+\Delta_S+D)\geq \kappa(K_S+\Delta_S+\epsilon D)=
\kappa\bigl( K_{S^m}+\Delta^m+\epsilon D^m\bigr)\geq 1.
\qed
$
%\medskip

\begin{lem}   \label{nef.vertucal.div.cor}
Let $g:X\to B$ be a proper morphism with connected general
fiber $X_g$.
Let $D$ be an effective, nef,  $\q$-Cartier $\q$-divisor on $X$.
 Then 
\begin{enumerate}
\item either $D|_{X_g}$ is a nonzero nef divisor
\item or $D$ is disjoint from $X_g$ and 
 $\bigl(D^{\dim B+1}\bigr)=0$. Thus $\nu(D)\leq \dim B$.
\end{enumerate}
\end{lem}

Proof. We are done if 
$D|_{X_g}$ is nonzero.
If it is zero then $D$ is vertical 
hence there is an ample divisor $L_B$ such that
$g^*L_B\sim D+E$ where $E$ is effective.  
Then
$$
\bigl(g^*L^r_B\bigr)-(D^r)=
\tsum_{i=0}^{r-1}\bigl(E\cdot g^*L^i_B\cdot D^{r-1-i}\bigr)
$$
shows that $  (D^r)  \leq  \bigl(g^*L^r_B\bigr) $.
Since $\bigl((g^*L_B)^{\dim B+1}\bigr)=g^*\bigl(L^{\dim B+1}_B\bigr)=0$, 
we conclude that
 $\bigl(D^{\dim B+1}\bigr)=0$.
 \qed

\begin{lem}  \label{nef.vertucal.div.lem}
Let $g:X\to B$ be a proper morphism with connected fibers
and $D$ an effective, exceptional,   $\q$-Cartier $\q$-divisor on $X$.
Then $D$ is not $g$-nef.
\end{lem}

Proof.  Let $|H|$ be a very ample
linear system on $X$ and $S\subset X$ the intersection of
$\dim X-2$ general members of $|H|$. Then
$g|_S:S\to B$ is generically finite over its image and
$D\cap S$ is   $g|_S$-exceptional. By the Hodge index theorem
we conclude that $\bigl(D^2\cdot H^{\dim X-2}\bigr)<0$, a contradiction.
\qed

\section{Asymptotic estimates for  cohomology groups}

\begin{say}\label{kodaira.formula}
 Let $X$ be a smooth variety and $g:X\to B$  a Calabi--Yau fiber space
of relative dimension $m$ over a smooth curve $B$.
Assume that $K_{X_g}\sim 0$ where $X_g$ denotes a general fiber. 
It is easy to see that the sheaves
$R^mg_*\o_X$ and $g_*\omega_{X/B}$
are line bundles and dual to each other. 
For elliptic surfaces these sheaves were computed by Kodaira.
His results were clarified and extended to higher dimensions by
\cite{fuj}. We will need the following consequences of their results.

The degree $\deg  g_*\omega_{X/B}$ is $\geq 0$ and it can be written as a 
sum of 2 terms. One is a global term
(determined by the $j$-invariant of the fibers in the elliptic case)
 which is zero iff  $g:X\to B$ is 
 {\it generically  isotrivial}, that is,
$g$ is  an analytically
 locally trivial fiber bundle
over a dense open set $B^0\subset B$. The other is a local term,
supported at the points where the local monodromy of the
local system $R^mg_*\q_{X^0}$ is nontrivial. There is a precise formula
for the local term, but we only need to understand 
what happens with generically isotrivial families. 
For these the local term is positive iff the local monodromy has
eigenvalue $\neq 1$ on  $g_*\omega_{X^0/B^0}\subset 
\o_{B^0}\otimes_{\q} R^mg_*\q_{X^0}$.

Over higher dimensional bases, $R^mg_*\o_X $ and $g_*\omega_{X/B} $
are rank 1 sheaves, and the above considerations describe their 
codimension 1 behavior. 
In particular, we see the following.
\begin{enumerate}
\item $c_1\bigl(g_*\omega_{X/B}\bigr)$ is linearly equivalent to
 a sum of effective $\q$-divisors. It is zero only if
 $g:X\to B$ is 
isotrivial over a dense open set $B^0$ and
the local   monodromy around each irreducible component of $B\setminus B^0$
 has
eigenvalue $= 1$ on  $g_*\omega_{X^0/B^0}\subset 
\o_{B^0}\otimes_{\q} R^mg_*\q_{X^0}$.
\item $c_1\bigl(R^mg_*\o_X\bigr)=-c_1\bigl(g_*\omega_{X/B}\bigr)$.
\end{enumerate}
Frequently $c_1\bigl(g_*\omega_{X/B}\bigr)$ is denoted by
$\Delta_{X/B}$.

\end{say}

\begin{cor} \label{main.chi.cor}
Let $g:X\to B$ be an elliptic fiber space
of dimension $n$
and $L$ a line bundle on $B$. Then
$$
\begin{array}{rcl}
\chi\bigl(X, g^*L^m\bigr)&=&
\frac{(L^{n-2}\cdot \Delta_{X/B})}{2(n-2)!} m^{n-2}+ O(m^{n-3}) \qtq{and}\\[1ex]
h^i \bigl(X, g^*L^m\bigr)&=&O(m^{n-3}) \qtq{for $i\geq 2$.}\\[1ex]
%h^0 \bigl(X, g^*L^m\bigr)&\geq &
%\frac{(L^{n-2}\cdot (\Delta+J)}{2(n-2)!} m^{n-2}+ O(m^{n-3}).
\end{array}
$$
\end{cor}

Proof. By the Leray spectral sequence,
$$
\chi\bigl(X, g^*L^m\bigr)=\sum (-1)^i\chi\bigl(B, L^m\otimes R^ig_*\o_X\bigr).
$$
For $i\geq 2$ the support of $R^ig_*\o_X $ has codimension $\geq 2$ in $B$,
hence its cohomologies
 contribute only to the $O(m^{n-3}) $ term.

Since $g$ has connected fibers, $g_*\o_X\cong \o_B$
and $c_1\bigl(R^1g_*\o_X\bigr)\simq -\Delta_{X/B}$ by %Kodaira's formula
(\ref{kodaira.formula}.2). 
We conclude by applying (\ref{2.term.RR.lem}) to both terms.\qed
\medskip

\begin{say}\label{22.say}
 Similar formulas apply to arbitrary Calabi--Yau fiber spaces
$g:X\to B$ with general fiber $F$.
For $m\gg 1$ we have
$$
H^i\bigl(X, g^*L^m\bigr)=H^0\bigl(B, L^m\otimes R^ig_*\o_X\bigr)=
\chi \bigl(B, L^m\otimes R^ig_*\o_X\bigr).
\eqno{(\ref{22.say}.1)}
$$
Setting $k=\dim B$, (\ref{2.term.RR.lem}) computes $H^i\bigl(X, g^*L^m\bigr) $
 as 
$$
\frac{m^{k}}{k!} h^i(F, \o_F)(L^{k})+
 \frac{m^{k-1}}{(k-1)!}
\Bigl(L^{k-1}\cdot \bigl(c_1(R^ig_*\o_X)-\tfrac{h^i(F, \o_F)}2 K_B\bigr)\Bigr)+
O(m^{k-2}).
$$
These imply that 
$$
\chi\bigl(X, g^*L^m\bigr)={\chi(F, \o_F)}\cdot \frac{m^{k}}{k!} (L^{k})+
O(m^{k-1}).
\eqno{(\ref{22.say}.2)}
$$
If $\chi(F, \o_F)\neq 0$ then this describes the asymptotic behavior  of
$\chi\bigl(X, g^*L^m\bigr) $.
However, if $\chi(F, \o_F)=0$, which happens for Abelian fibers,
then we have to look at the next term which gives that
$$
\chi\bigl(X, g^*L^m\bigr)=
 \frac{m^{k-1}}{(k-1)!}
\Bigl(L^{k-1}\cdot \tsum_{i=1}^{\dim F} (-1)^i c_1(R^ig_*\o_X)\Bigr)+
O(m^{k-2}).
\eqno{(\ref{22.say}.3)}
$$
If $F$ is an elliptic curve then the sum on the right hand side
has only 1 nonzero term. For higher dimensional Abelian fibers
there are usually several nonzero terms and sometimes they cancel each other.

This is one reason why  elliptic fibers are easier to study than
higher dimensional Abelian fibers. The other difficulty 
with higher dimensional  fibers is that the Euler characteristic only tells us
that  $h^0+h^2+h^4+\cdots$ grows as expected.
Proving that $h^0\neq 0$  would need additional arguments.
\end{say}

The next result, while stated in all dimensions,
is truly equivalent to Kodaira's formula \cite[V.12.2]{bpv}.

\begin{cor} \label{td2.cor}
Let $g:X\to B$ be a relatively minimal  elliptic fiber space
of dimension $n$
and $L$ a line bundle on $B$.
Then
$\bigl(L^{n-2}\cdot \Delta_{X/B}\bigr)=
\bigl(g^*L^{n-2}\cdot \operatorname{td}_2(X)\bigr)$.
\end{cor}

Proof. Expanding the Riemann--Roch formula
$\chi(X,L)=\int_X \operatorname{ch}(L)\cdot \operatorname{td}(X)
$
gives that 
$$
\chi\bigl(X, L^m\bigr)=
\frac{(L^n)}{n!}\cdot m^n-
\frac{(L^{n-1}\cdot K_X)}{2(n-1)!}\cdot  m^{n-1}+ 
\frac{(L^{n-2}\cdot \operatorname{td}_2(X))}{(n-2)!}\cdot m^{n-2}+ O(m^{n-3}).
$$
Comparing this with (\ref{main.chi.cor})
yields the claim.\qed
\medskip

We used several versions of the asymptotic
Riemann--Roch formula.

\begin{lem} \label{2.term.RR.lem}
Let $X$ be a normal, projective variety 
of dimension $n$, $L$ a line bundle
on $X$ and $F$ a coherent sheaf that is locally free
in codimension 1. Then
$$
\chi\bigl(X, L^m\otimes F\bigr)=
\frac{(L^n)\cdot\rank F}{n!} m^n+
\frac{\bigl(L^{n-1}\cdot (c_1(F)-\frac{\rank F}2 K_X)\bigr)}{(n-1)!} m^{n-1}+
 O(m^{n-2}).\qed
$$
\end{lem}

\begin{say}[Riemann--Roch with rational singularities]\label{RR.ratsing.say}
The Todd classes of a singular variety $X$ are not always easy to
compute, but if $X$ has  rational singularities
then there is a straightforward formula in terms of the
Chern classes of any resolution  $h:X'\to X$.

By definition, rational singularity
means that $R^ih_*\o_{X'}=0$ for $i>0$. Thus
$\chi(X,L)=\chi(X', h^*L)$  for any line bundle
$L$ on $X$. By the projection formula this implies that
$\chi(X,L)=\int_X \operatorname{ch}(L)\cdot h_*\operatorname{td}(X')
$
and in fact $ \operatorname{td}(X)= h_*\operatorname{td}(X')$
(cf.\ \cite[Thm.18.2]{MR1644323}.)
In particular we see that
 the second Todd class of $X$ is
 $$
\operatorname{td}_2(X)=h_* \Bigl(\frac{c_1(X')^2+c_2(X')}{12}\Bigr).
$$
\end{say}

The following numerical version of (\ref{main.chi.cor})
was used in the proof of  (\ref{main.char.thm.td2}).

\begin{lem} \label{there.are.sections.lem}
Let $X$ be a normal, projective variety of dimension $n$.
Let  $L$  be a nef line bundle on $X$ such that
$(L^n)=(L^{n-1}\cdot K_X)=0$ but $(L^{n-1} )\neq 0$.
Then
$$
h^0(X, L^m)-h^1(X, L^m)=
\frac{(L^{n-2}\cdot \operatorname{td}_2(X))}{(n-2)!}\cdot m^{n-2}+O(m^{n-3}).
$$
\end{lem}

Proof. The assumptions $(L^n)=(L^{n-1}\cdot K_X)=0$ imply
that the right hand side equals $\chi(X, L^m)$. 
Thus the equality follows if 
$h^i(X, L^m)=O(m^{n-3})$ for $i\geq 2$.
The latter is a special case of (\ref{growth.of.coh.groups}).\qed

\begin{lem} \label{growth.of.coh.groups}
Let $X$ be a projective variety  of dimension $n$
and $F$ a torsion free
coherent sheaf on $X$. Let $L$ be a nef line bundle on $X$
and set $d=\nu(X,L)$. Then
$$
\begin{array}{lcl}
h^i(X, F\otimes L^m)&=& O(m^d)\qtq{for $i=0,\dots, n-d$ and}\\
h^{n-j}(X, F\otimes L^m)&=& O(m^{j-1}) \qtq{for $j=0,\dots, d-1$.}
\end{array}
$$
\end{lem}

Note the key feature of the estimate: the order of growth of $H^i$ is
$m^d$ for $i\leq n-d$, then for $i=n-d+1$ it drops by 2 to
$m^{d-2}$ and then it drops by 1 for each increase of $i$.
This strengthens \cite[1.4.40]{laz-book} but the proof
is essentially the same.
\medskip

%The estimates are equivalent to  NO!!!!!!!!!  i=n-d
%$$
%h^i(X, F\otimes L^m)= O\bigl(m^{\min\{d, n-i-1\}}\bigr) \qtq{for all $i$.}
%$$

%If there is an exact sequence $0\to F_1\to F_2\to F_3\to 0$
%and the lemma holds for $F_1, F_3$ then it clearly holds for $F_2$. 
%Thus it is sufficient to consider the case when  $F$ is torsion free. 
%We may even assume that 
%$F\cong \o_X(-D)$ for some effective Cartier divisor $D$.

Proof. We use  induction on $\dim X$.
By Fujita's theorem (cf.\ \cite[1.4.35]{laz-book})
 we can choose a general very ample 
divisor $A$ on $X$ such that
 $$
h^i(X, F\otimes \o_X(A)\otimes L^m)= 0 \qtq{for all $i\geq 1$ and $m\geq 1$.}
$$
We get an exact sequence
$$
0\to F\otimes L^m \to F\otimes \o_X(A)\otimes L^m \to
G\otimes L^m \to 0
$$
where $G$ is  a torsion free
coherent sheaf  on $A$. For $i\geq 1$ its long cohomology sequence
gives surjections (even isomorphisms for $i\geq 2$)
$$
H^{i-1}(A, G\otimes L^m )\onto H^i(X, F\otimes L^m).
$$
By induction this shows the claim except for $i=0$. 

One can realize $F$ as a subsheaf of a sum of line bundles,
thus it remains to prove that 
 $H^0(X, F\otimes L^m)=O(m^d)$ 
when $F\cong \o_X(H)$ is a very ample line bundle. 
The exact sequence
$$
0\to L^m\to \o_X(H)\otimes L^m\to \o_H(H|_H)\otimes L^m \to 0
$$
finally reduces the problem to 
$\kappa(L)\leq \nu(L)$ which was discussed in (\ref{kappa.nu.defn}).
 \qed

\section{Deforming morphisms}

Here we answer Question \ref{ques1} but first two technical issues
need to be discussed: the distiction between
\'etale and quasi-\'etale covers and the
existence of non-Calabi--Yau deformations.
Both appear  only for singular
  Calabi--Yau varieties.

\begin{defn} Following \cite{cat-qed}, a finite morphism
$\pi:U\to V$ is called {\it quasi-\'etale} if there is a 
closed subvariety $Z\subset V$ of codimension $\geq 2$
such that $\pi$ is  \'etale over $V\setminus Z$.

If $V$ is a normal variety, then there is a one-to-one
correspondence between quasi-\'etale covers of $V$
and finite, \'etale covers of $V\setminus \sing V$.

In particular,
if $X$ is a Calabi--Yau variety then there is a 
quasi-\'etale morphism $ X_1\to X$ such that
$K_{X_1}\sim 0$. 

Among all such covers $X_1\to X$ there is a unique smallest one,
called the {\it index 1 cover} of $X$, which is Galois with cyclic Galois group.
We  denote it by $X^{\rm ind}\to X$. 
\end{defn}

\begin{say}[Deformation theory] For a general introduction, see
\cite{har-def}. By a deformation of a proper scheme (or analytic space)
$X$ we mean
a flat, proper  morphism  $g:{\mathbf X}\to (0\in S)$
to a pointed scheme (or analytic space) together with a fixed
isomorphism $X_0\cong X$. 

By a deformation of a morphism of proper schemes (or analytic spaces)
$f:X\to Y$ we mean a morphism 
${\mathbf f}:{\mathbf X}\to {\mathbf Y}$
where ${\mathbf X} $ is a deformation of $X$,
${\mathbf Y} $ is a deformation of $Y$
and ${\mathbf f}|_{X_0}=f$. 

When we say that an assertion holds for all
{\it small deformations} of $X$, this means that for every
deformation $g:{\mathbf X}\to (0\in S)$ there is an \'etale
(or analytic) neighborhood  $(0\in S')\to (0\in S)$ such that
the assertion holds for  
$g':{\mathbf X}\times_SS'\to (0\in S')$.
\end{say}

\begin{say}[Deformations of Calabi--Yau varieties]\label{CY.def.say.say}
Let $X$ be a Calabi--Yau variety. If $X$ is smooth
(or has canonical singularities, or $K_X$ is Cartier)
then every small defomation of
$X$ is again a Calabi--Yau variety. This, however, fails in general;
see (\ref{non.CY.def.exmp})
 for an example where $X$ is a surface with quotient singularities.

Dealing with such unexpected deformations is a basic problem
in the moduli theory of higher dimensional varieties;
see \cite[Sec.4]{k-modsurv}, \cite[Sec.14B]{hac-kov} or \cite{abr-hass}
for a discussion and solutions.
For Calabi--Yau varieties one can use a global trivialization of
the canonical bundle to get a much simpler answer.

We say that a defomation $g:\mathbf X\to (0\in S)$ of $X$ 
over a reduced, local  space $S$ 
is a {\it Calabi--Yau deformation} if 
the following equivalent conditions hold:
\begin{enumerate}
\item Every fiber of $g$ is a Calabi--Yau variety.
\item The deformation can be lifted to a deformation  
$g^{\rm ind}:{\mathbf X}^{\rm ind}\to (0\in S)$  of $X^{\rm ind}$, 
the  index 1 cover of $X$.
\end{enumerate}

% The following variant is much easier to use in practice.
% We adopt it as our definition.
% Let $(0\in S)$ be a local scheme. A flat defomation
% $g:{\mathbf X}\to (0\in S)$  of $X$ is a {\it Calabi--Yau deformation} if
% it can be lifted to a deformation  
% $g^{\rm ind}:{\mathbf X}^{\rm ind}\to (0\in S)$  of $X^{\rm ind}$, 
% the  index 1 cover of $X$.

Thus studying  Calabi--Yau deformations of Calabi--Yau varieties
is equivalent to  studying deformations of Calabi--Yau varieties
whose canonical class is Cartier. As we noted, for the latter
every deformation is automatically a  Calabi--Yau deformation.
Thus we do not have to deal with this issue at all.
\end{say}

%  Not every deformation of an
% elliptic  Calabi--Yau variety is elliptic.
% For instance,  a general deformation of 
%  $(\mbox{Abelian variety})\times(\mbox{elliptic curve})$
%  has no elliptic fiber space structure.
% Similarly, every elliptic K3 surface has non-elliptic deformations.
% The next result says that these are essentially the only
% such examples.

\begin{thm}\label{smooth.horik.thm}
Let $X$ be  a  Calabi--Yau variety and
$g:X\to B$  an elliptic fiber space. 
 Then at least one of the following holds.
\begin{enumerate}
\item The morphism $g$ extends to every  small  Calabi--Yau deformation of $X$.
\item There is a quasi-\'etale cover $\tilde X\to X$ such that
the Stein factorization $\tilde g:\tilde X\to \tilde B$
of $\tilde X\to B$ is one of the following
\begin{enumerate}
\item  $\bigl(\tilde g:\tilde X\to \tilde B\bigr)\cong
\bigl(p_1: \tilde B\times(\mbox{elliptic curve})\to \tilde B\bigr)$ 
where $p_1$ is the first projection or
\item  $\bigl(\tilde g:\tilde X\to \tilde B\bigr)\cong
\bigl(p_1: \tilde Z\times(\mbox{elliptic K3})\to \tilde Z\times \p^1\bigr)$
\end{enumerate}
where $\tilde Z$ is a Calabi--Yau variety of dimension $\dim X-2$
and $p_1$ is the product of the first projection with the
elliptic pencil map of the K3 surface.
\end{enumerate}
\end{thm}

 Proof. As noted in (\ref{CY.def.say.say}), we may assume that $K_X\sim 0$.

By \cite{KMM92a} there is a unique map
(up-to birational equivalence)
 $h:B\map Z$  whose general fiber $F$ is smooth, proper, rationally connected
and whose target $Z$ is not uniruled by \cite{ghs}.
(See \cite[Chap.IV]{rc-book} for a detailed treatment
or \cite{MR2011743} for an introduction.)
Next apply \cite[Thm.3]{k-lars} to $X\map Z$
to conclude that there is a finite \'etale cover $\tilde X\to X$, 
 a product decomposition   $\tilde X\cong  Y\times \tilde Z$
and a generically finite map $\tilde Z\map Z$  that factors
 $\tilde X\map Z$.

If $\dim Z=\dim B$ then we are in case (2.a).
If $\dim Z=\dim B-1$ then the generic fiber of $\tilde B\to \tilde Z$
is $\p^1$. Furthermore, $\dim Y=2$, hence either
$Y$ is an elliptic  K3 surface
 and we are in case (2.b)
or $Y$ is an Abelian surface that has an elliptic pencil and
after a further cover we are again in case (2.a).

It remains to prove that if $\dim F\geq 2$ then 
the assertion of (1) holds.
By (\ref{horik.thm})
it is sufficient to check that 
$$
\Hom_B\bigl(\Omega_B,R^1g_*\o_X\bigr)=0.
$$
To see this, note that   (\ref{kodaira.formula}) and $\omega_X\sim \o_X$
imply that 
$R^1g_*\o_X\cong \bigl(g_*\omega_{X/B}\bigr)^{-1}\cong \omega_B$.
Therefore
$$
\shom_B\bigl(\Omega_B,R^1g_*\o_X\bigr)\cong 
\shom_B\bigl(\Omega_B,\omega_B\bigr)
\cong \bigl(\Omega_B^{\dim B-1}\bigr)^{**}
$$
where $({\ })^{**}$ denotes the double dual or reflexive hull.
By taking global sections we get that
$$
\Hom_B\bigl(\Omega_B,R^1g_*\o_X\bigr)=
H^0\bigl(B, \bigl(\Omega_B^{\dim B-1}\bigr)^{**}\bigr).
$$
Let $B'\to B$ be a resolution of singularities and
$F'\subset B'$ a general fiber of $B'\map Z$. 
Since $F'$ is rationally connected, if $C\subset F'$ is a general rational curve
then  
$$
T_{F'}|_C\cong \tsum \o_C(a_i)\qtq{where $a_i>0$ $\forall\, i$;}
$$
see \cite[IV.3.9]{rc-book}. Thus
$T_{B'}|_C$ is a sum of line bundles $\o_C(a_i)$
where $a_i>0$ for $\dim F$ summands and $a_i=0$ for the rest.
Since $\dim F\geq 2$ we conclude that
$$
\wedge^{\dim B-1}T_{B'}|_C\cong \tsum \o_C(b_i)\qtq{where $b_i>0$ for every  $i$.}
$$
By duality this gives that
$H^0\bigl(B', \Omega_{B'}^{\dim B-1}\bigr)=0$.
Finally we use that $B$ has log terminal singularities by \cite{naka88}
and so \cite{GKKP10} shows that 
$$
\Hom_B\bigl(\Omega_B,R^1g_*\o_X\bigr)=
H^0\bigl(B, \bigl(\Omega_B^{\dim B-1}\bigr)^{**}\bigr)=
H^0\bigl(B', \Omega_{B'}^{\dim B-1}\bigr)=0.\qed
$$

% Thus $H^0\bigl(B, T_B\otimes R^1g_*\o_X\bigr)=0$

% Note that   (\ref{kodaira.formula}) and $\omega_X\sim \o_X$
% imply that 
% $R^1g_*\o_X\cong \bigl(g_*\omega_{X/B}\bigr)^{-1}\cong \omega_B.$
% Therefore
% $$
% T_B\otimes R^1g_*\o_X\cong T_B\otimes \omega_B\cong \Omega_B^{\dim B-1}.
% $$
% Thus $H^0\bigl(B, T_B\otimes R^1g_*\o_X\bigr)=0$
% and so (\ref{horik.thm})
% implies that every  small deformation of $X$ is also  an elliptic fiber space.
% \qed

%\medskip

We are now ready to answer Question \ref{ques1}.

\begin{thm}\label{thm.ques1.general} Let $X$ be an elliptic Calabi--Yau
variety such that $H^2(X, \o_X)=0$. Then 
every small Calabi--Yau deformation of $X$ is also an 
elliptic Calabi--Yau variety.
\end{thm}

Proof. Let $g:X\to B$ be an  elliptic Calabi--Yau
variety.
By (\ref{smooth.horik.thm}) every small Calabi--Yau deformation of $X$
 is also an 
elliptic Calabi--Yau variety except possibly when
there is a quasi-\'etale cover $\tilde X\to X$ such that
\begin{enumerate}
\item either $\tilde X\cong \tilde Z\times (\mbox{elliptic curve})$
\item or $\tilde X\cong \tilde Z\times (\mbox{elliptic K3})$.
\end{enumerate}
In both cases, $\tilde X$ can have non-elliptic deformations
but we show that these do not correspond to a deformation of $X$.
Here we use that $H^2(X, \o_X)=0$.

Let $\pi:{\mathbf X}\to (0\in S)$ be a flat deformation of $X$
over a local scheme $S$.
Let $L$ be the pull-back of an ample line bundle from $B$ to $X$.
Since $H^2(X, \o_X)=0$, $L$ lifts to a line bundle
${\mathbf L} $ on $ {\mathbf X}$   (cf.\ \cite[p.236-16]{FGA})
thus we get a
line bundle $\tilde {\mathbf L} $ on $ \tilde {\mathbf X}$.
We need to show that a large multiple of  $ {\mathbf L} $
is base-point-free over $S$; then it
gives the required morphism  
 ${\mathbf g}: {\mathbf X}\to {\mathbf B}$.
One can check base-point-freeness of some multiple after
a finite surjection, thus it is enough to show
that some multiple of  $\tilde {\mathbf L} $ 
is base-point-free over $S$.

The first case (more generally, deformations of products
with Abelian varieties) is treated in 
(\ref{ab.acts.deforms}). 

In the K3 case note first that every small deformation of
$\tilde X$ is of the form $\tilde {\mathbf Z}\times_S \tilde {\mathbf F}$
where $\tilde {\mathbf F}\to S$ is a flat family of K3 surfaces.
This is a trivial case of (\ref{horik.thm}); see
(\ref{deform.conj.say}) for an elementary argument. 
Hence we only need to show that the  restriction
of  $\tilde {\mathbf L} $ to $ \tilde {\mathbf F}$
is base-point-free over $S$. Equivalently, that the elliptic structure of the
central K3 surface  $\tilde F$ is preserved by our deformation.
 The restriction
of  $\tilde {\mathbf L} $ to every fiber of $\tilde {\mathbf F}\to S$
 gives a nonzero, nef line bundle
with self-intersection 0, hence an elliptic pencil.\qed
\medskip

\begin{say}[Deformation of sections]
Let $g:X\to B$ be an elliptic Calabi--Yau fiber space
with a section $S\subset X$. Let us assume first that $S$ is a Cartier
divisor in $X$. (This is automatic if $X$ is smooth.)
Then $S$ is $g$-nef,  $g$-big and $S\sim_{\q, g} K_X+S$ hence
$R^ig_*\o_X(S)=0$ for $i>0$; cf.\ \cite[Sec.2.5]{km-book}.
Thus  $H^i(X, \o_X(S))=H^i(B, g_*\o_X(S))$ for every $i$.
In order to compute $g_*\o_X(S)$ we use the exact sequence
$$
0\to \o_B=g_*\o_X\stackrel{\alpha}{\to} g_*\o_X(S) \to g_*\o_S(S|_S)
$$
A degree 1 line bundle over an elliptic curve has only 1 section, thus
 $\alpha$ is an isomorphism over an open set where
the fiber is a smooth elliptic curve. Since $g_*\o_S(S|_S)\cong \o_S(S|_S)$
is torsion free we conclude that $g_*\o_X(S)\cong \o_B$. Thus
$$
H^1(X, \o_X(S))=H^1(B, \o_B)\subset H^1(X, \o_X).
$$
If $H^2(X, \o_X)=0 $ then the line bundle $\o_X(S)$ lifts to
every small deformation of $X$ and if $H^1(X, \o_X)=0 $
then the unique section of $\o_X(S)$ also lifts.

The situation is quite different if the section is not assumed Cartier.
For instance, let $X_0\subset \p^2\times \p^2$ be a general hypersurface
of multidegree $(3,3)$ containing  $S:=\p^2\times \{p\}$ for some 
point $p$. Then $X_0$ is a Calabi--Yau variety and the first projection
 shows that it is elliptic with a section. Note that $X_0$ is singular,
it has 9 ordinary nodes along $S$.

By contrast, if $X_t\subset \p^2\times \p^2$ is a smooth hypersurface
of multidegree $(3,3)$ then the restriction map
$\pic\bigl(\p^2\times \p^2\bigr)\to \pic(X_t)$
is an isomorphism by the Lefschetz hyperplane theorem.
 Thus the degree of every divisor $D\subset X_t$ 
on the general fiber of the first projection $X_t\to \p^2$
is a multiple of 3. Therefore $X_t\to \p^2$ does not even have a 
rational section.
\end{say}

As an aside, we consider the general question of
deforming morphisms $g:X\to Y$ whose target
is not uniruled. 

There are some obvious examples
when not every deformation of $X$ gives a deformation of $g:X\to Y$.
For example, let $A_1, A_2$ be positive dimensional Abelian varieties
and  $g:A_1\times A_2 \to A_2$ the second projection.
A general deformation of $A_1\times A_2 $ 
is a simple Abelian variety which has no maps to
lower dimensional Abelian varieties. 
One can now get more complicated examples
by replacing $A_1\times A_2 $ by say a complete intersection
subvariety or by a cyclic cover.
The next result says that this essentially gives
all examples.

\begin{thm}\label{smooth.horik2.thm}
Let $X$ be  a projective variety with rational singularities,
$Y$ a normal variety
 and
$g:X\to Y$  a surjective morphism with connected fibers. 
Assume that $Y$ is not uniruled. Then at least one of the following holds.
\begin{enumerate}
\item Every  small deformation of $X$ gives a deformation of $(g:X\to Y)$.
\item There is a quasi-\'etale cover $\tilde Y\to Y$,
a smooth variety $Z$ and positive dimensional Abelian varieties
$A_1, A_2$  such that
the lifted morphism
 $\tilde g:\tilde X:=X\times_Y \tilde Y\to  \tilde Y $ factors as
$$
\begin{array}{ccc} 
\tilde X  &  \to & Z\times A_2\times A_1\\
\tilde g\downarrow\hphantom{\tilde g}  && \downarrow \\
\tilde Y & \cong & Z\times A_2.
\end{array}
$$

\end{enumerate}
\end{thm}

Proof.  By (\ref{horik.thm}) 
 every   deformation of $X$ gives a deformation of $g:X\to Y$ if
$$
\Hom_Y\bigl(\Omega_Y, R^1g_*\o_X\bigr)=0.
\eqno{(\ref{smooth.horik2.thm}.3)}
$$
Thus we need to show that if (\ref{smooth.horik2.thm}.3) fails
then we get a structural description as in (\ref{smooth.horik2.thm}.2).

Let $C\subset Y$ be a very general complete intersection curve.
Since $Y$ is not uniruled, $\Omega_Y|_C$ is semi-positive by
\cite{MR949837}; see also \cite[Sec.9]{k-etal}.

Let $Y^0\subset Y$ be a dense open set 
and $X^0:=g^{-1}(Y^0)$ 
such that
$g^0:X^0\to Y^0$ is smooth. Set $C^0:=X^0\cap C$. 
By \cite{Steenbrink75}, $\bigl(R^1g_*\o_X\bigr)|_C$
is the (lower) canonical extension of the top quotient of the
variation of Hodge structures $R^1g^0_*\q_{X^0}|_{C^0}$.
(Note that \cite{Steenbrink75} works with $\omega_{X^0/Y^0}$
but the proof is essentially the same; see
\cite[pp.177--179]{k-hdi2}.)
Thus $\bigl(R^1g_*\o_X\bigr)|_C$ is  semi-negative
by \cite{Steenbrink75}. 
Moreover, the  part that is not strictly negative
corresponds to a  variation of sub-Hodge structures
that  is a direct summand and becomes trivial
after a suitable quasi-\'etale cover $Y_1\to Y$ \cite[Thm.4.2.6]{MR0498551}.
This direct summand 
corresponds to a  direct factor of the Albanese variety of $X_1:=Y_1\times_YX$,
giving the Abelian variety $A_1$. 

Once the flat part of $R^1g^0_*\q_{X^0}$ is trivial,
$$
\shom_Y\bigl(\Omega_Y, R^1g_*\o_X\bigr)|_C\cong
\bigl(T_Y\otimes R^1g_*\o_X\bigr)|_C
$$
 has a global section iff $T_Y|_C$
has a global section. Thus $H^0(Y, T_Y)\neq 0$
and so $\dim \aut(Y)>0$. Since $Y$ is   not uniruled,
$\aut^0(Y)$ has no linear algebraic subgroups, thus 
the connected component $\aut^0(Y)$ is an Abelian variety $A_2$.
By (\ref{ab.action.et.triv}), $A_2$ becomes a direct factor
after a suitable \'etale cover $\tilde Y\to Y_1\to Y$. \qed 
\medskip

The following result was essentially known to
\cite{serre-1958, MR0164973}; see
\cite{brion-2007} for the general theory.

\begin{prop}\label{ab.action.et.triv}
 Let $W$ be a normal, projective variety
and $A$ an Abelian variety acting faithfully on $W$. Then there is a
 normal, projective variety $Z$ and an $A$-equivariant
\'etale morphism $A\times Z\to W$.
\end{prop}

Proof. Let $T\subset W$ be the generic orbit.
The quotient $V:=W/A$ exists; it is the normalization of the closure of
the point $[T]$ corresponding to $T$ in the Chow variety $ \chow(W)$.
Since $W\to V$ is a generically isotrivial $A$-bundle,
 using (\ref{isortiv.standard.form.thm}),
we obtain that 
$W\cong \bigl(\tilde V\times A\bigr)/G$ 
where $\tilde V\to V$ is a  finite (ramified cover)
and $G$ acts faithfully on $\tilde V$ and on $A$. 
Since the $A$-action descends to $W$, the $G$-action on $A$
commutes with translations. 
An automorphism of an Abelian variety that commutes with
all translations is itself a translation.
Thus $G$ acts on $A$ via translations
and so  the $G$-action on $\tilde V\times A$
is fixed point free. Therefore
$\tilde V\times A\to W$ is \'etale.\qed

\medskip

The following is a combination of \cite[Thm.8.1]{MR0417458}
and the method of \cite[Thm.8.2]{MR0417458} in the smooth case
and \cite[Prop.3.10]{bhps} in general.

\begin{thm}\label{horik.thm}
Let $f:X\to Y$ be a morphism of proper varieties
such that $f_*\o_X=\o_Y$ and 
$\Hom_Y\bigl(\Omega_Y,R^1f_*\o_X\bigr)=0$.

Then for  every small deformation of $\mathbf X$ of $X$ 
there is a small deformation  $\mathbf Y$ of $Y$
such that $f$ lifts to
${\mathbf f}:{\mathbf X}\to {\mathbf Y}$.\qed
\end{thm}

Note that  if $X$ is smooth or, more generally, it it
 has rational singularities, then 
$R^1f_*\o_X $ is a reflexive sheaf by \cite[7.8]{k-hdi1}. Thus,
if $Y$ is normal then
$$
\Hom_Y\bigl(\Omega_Y,R^1f_*\o_X\bigr)=
\Hom_{Y^{\rm ns}}\bigl(\Omega_{Y^{\rm ns}},R^1f_*\o_X|_{Y^{\rm ns}}\bigr)
$$
where $Y^{\rm ns}\subset Y$ is the smooth locus.
Therefore we can check the vanishing of 
$\Hom_Y\bigl(\Omega_Y,R^1f_*\o_X\bigr)=0$ by finding general
projective curves $C\subset Y^{\rm ns}$ such that the vector bundle
$\bigl(T_{Y}\otimes R^1f_*\o_X\bigr)|_C$ has no global sections.

\section{Smoothings of very singular varieties}

One can frequently construct smooth varieties by first
exhibiting some very singular, even reducible schemes
with suitable numerical invariants and then smoothing them.
For such Calabi--Yau examples see \cite{MR1296351}.
Thus it is of interest to know when an elliptic fiber space
structure is preserved by a smoothing.
In some cases, when (\ref{thm.ques1.general}) does not apply,  the
following result,  relying on 
(\ref{main.char.thm.cor.q2}),  provides a quite satisfactory answer.

\begin{prop}\label{def.sing.say.prop}
Let $X$ be a  projective, Gorenstein  scheme of pure dimension $n$
such that $\omega_X$ is  numerically trivial and  $H^2(X, \o_X)=0$.
Let $g:X\to B$ be a morphism whose general fibers 
(over every irreducible component of $B$)
are curves of arithmetic genus 1.
Assume also that  every irreducible component of $X$
dominates an irreducible component of $B$.

Let $L_B$ be an ample  line bundle on $B$ 
and assume that $\chi\bigl(X, g^*L^m_B\bigr)$
is a polynomial of degree  $\dim X-2$.
Then  every  smoothing of $X$ is  an   elliptic fiber space. 
\end{prop}

{\it Warning.} 
Note that we do not claim that $g$ lifts to every deformation of $X$.
In the example (\ref{rtle.elliptic.bad.exmp})
$X$ has smoothings, which are elliptic, and also other singular deformations
that are not elliptic.
\medskip

Proof.  As before, $H^2(X, \o_X)=0$ implies that 
$g^*L_B $ lifts to every small deformation \cite[p.236-16]{FGA}. 
Thus we have a 
deformation  $h:\bigl({\mathbf X}, {\mathbf L}\bigr) \to (0\in S)$ 
of $(X_0,L_0)\cong(X,L=g^*L_B)$.

We claim that ${\mathbf L} $  is $h$-nef and 
$K_{\mathbf X}$ is trivial on the fibers of $h$.
This is a somewhat delicate point since
 being nef is not known to be an
open condition in general. We go around this problem as follows.

Let $\bigl(X_{gen}, L_{gen}\bigr)$ be a generic fiber.
(Note the difference between generic and general.)
First we show  that  $L_{gen} $ is nef and $K_{X_{gen}}\simq 0$.
 Indeed, assume that 
$\bigl(L_{gen}\cdot C_{gen}\bigr)<0$ for some curve $C_{gen}$.
Let $C_0\subset X_0$ be a specialization of  $C_{gen}$. Then
$\bigl(L_0\cdot C_0\bigr)=\bigl(L_{gen}\cdot C_{gen}\bigr)<0$
gives a contradiction. A similar argument shows that
$\bigl(K_{X_{gen}}\cdot C_{gen}\bigr)=0$ for every curve $C_{gen}$.

Next,  the deformation invariance of $\chi\bigl(X, g^*L^m_B\bigr)$ and 
Riemann--Roch (cf.\ (\ref{td2.cor}) and (\ref{RR.ratsing.say})) show that
$$
\bigl(L_{gen}^{n-2}\cdot c_2(X_{gen})\bigr)=(n-2)!\cdot 
(\mbox{coefficient of $m^{n-2}$ in $\chi\bigl(X, g^*L^m_B\bigr)$}).
$$
Therefore $\bigl(L_{gen}^{n-2}\cdot c_2(X_{gen})\bigr)>0$
and, as we noted after  (\ref{main.char.thm.td2}), this implies that 
 $|mL_{gen}|$ is base point free for some $m>0$.

Thus  there is a dense Zariski open subset $S^0\subset S$ such that 
$|mL_s|$ is base point free for $s\in S^0$, hence 
$\bigl(X_s, L_s\bigr)$ is  an elliptic fiber space
for $s\in S^0$. 
We repeat the argument for the generic points of $S\setminus S^0$ and 
conclude by Noetherian induction. 
\qed
\medskip

Note that the universal deformation space of a proper scheme
can be represented by a scheme \cite{Artin69b}, 
thus the above argument takes care of analytic deformations
as well. It may be useful, however, to see how to modify
the proof to work directly in the analytic case when 
there are no generic points.

The (Barlet or Douady) space of curves in 
$h:{\mathbf X} \to (0\in S)$ has only countably many
irreducible components, thus  
 there are countably many
closed subspaces  $S_i\subsetneq S$ such that
every curve $C_s\subset X_s$ is deformation equivalent to a curve
$C_0\subset X_0$. In particular,
$L_s$ is nef and $K_{X_s}\simq 0$ whenever $s\not\in \cup S_i$. 
Thus $\bigl(X_s, L_s\bigr)$ is  an elliptic fiber space
for  $s\not\in \cup S_i$.

By semicontinuity, there are closed subvarieties  $T_m\subsetneq S$ such that
$$
h_*\o_{\mathbf X}(m{\mathbf L})\otimes \c_s=H^0\bigl(X_s, \o_{X_s}(mL_s)\bigr)
\qtq{for $s\not\in T_m$. }
$$
 Thus if $s\not\in \cup_i S_i\bigcup \cup_m T_m$
and $ \o_{X_s}(m_0L_s)$ is generated by global sections then
$$
\phi_{m_0}: h^*\bigl(h_*\o_{\mathbf X}(m_0{\mathbf L})\bigr)\to 
\o_{\mathbf X}(m_0{\mathbf L})
$$
is surjective along $X_s$. Thus 
there is a dense Zariski open subset $S^0\subset S$ such that 
$ \phi_{m_0}$ is surjective for all $s\in S^0$.
Now we can finish  by Noetherian induction as before.

% The proof of
% (\ref{main.char.thm.cor.q1}) shows that, over the generic point,
% $L_{gen}$ is nef and $K_{X_{gen}}$ is numerically trivial.

% and we finish as in the  proof of
% (\ref{main.char.thm.cor.q1}).\qed
% \medskip

% Proof. The condition $H^2(X, \o_X)=0$ guarantees that
% $L$ lifts to  every small deformation of $X$  \cite[p.236-16]{FGA}.
% All the assumptions of (\ref{main.char.thm.cor.q2}) are invariant under
% smooth deformations of a pair $(X,L)$, with the exception of
% $L$ being nef, which is not known to be an
% open condition in general. We go around this problem as follows.

% We start with algebraic deformations.

% For analytic deformations the argument can be modified as follows.

% \begin{cor}\label{main.char.thm.cor.q1}
%   Let $X$ be a smooth, projective  variety of dimension $n$
% and $L$ a Cartier divisor on $X$ such that
% $(X,L)$ is an  elliptic fiber space. 
% Assume that  $H^2(X, \o_X)=0$, $K_X\simq 0$ and 
% $\bigl(L^{n-2}\cdot c_2(X)\bigr)>0$.
% Then  every small deformation of $X$ is  also an   elliptic fiber space.
% \end{cor}

\section{Calabi--Yau orbibundles}

The techniques  of  this section are mostly taken from
\cite[Sec.6]{k-shaf} and  \cite{nakayama}.

\begin{defn}\label{isortiv.standard.form.exmp}
A Calabi--Yau fiber space  $g:X\to B$ is called an
{\it orbibundle} if it can be obtained by the
following construction.
 
 Let $\tilde B$ be a normal variety, $F$ a Calabi--Yau variety
and $\tilde X:=\tilde B\times F$.
Let $G$ be a finite group,  $\rho_B:G\to \aut(\tilde B)$
and   $\rho_F:G\to \aut(F)$ two faithful representations.
Set
$$
\bigl(g:X\to B\bigr):=
\bigl(\tilde X/G {\to} \tilde B/G\bigr);
$$
it is a generically isotrivial
Calabi--Yau fiber space.

(It would seem more natural to require the above
property only locally on $B$. We see in (\ref{isortiv.standard.form.thm})
that in the algebraic case the two version are equivalent.
However, if $X$ is a K\"ahler manifold, then the local
and global versions are different.)

For any non-ruled variety $Z$, the connected component $\aut^0(Z)$ of
$\aut(Z)$ is an Abelian variety, its elements are called translations.
The quotient $\aut(Z)/\aut^0(Z)$ is the discrete part of the
automorphism group.

For $G$ acting on $F$, 
let $G_t:=\rho_F^{-1}\aut^0(F)\subset G$ be the 
normal subgroup of translations
 and set $X^d:=\tilde X/G_t$. 
Then  $G_d:=G/G_t$ acts on $X^d$ and $X=X^d/G_d$. 
Thus every orbibundle comes with 2 covers:
$$
\begin{array}{ccccc}
X & \stackrel{\tau_X}{\longleftarrow} & X^d  & \stackrel{\pi_X}{\longleftarrow} & \tilde X\\
g\downarrow\hphantom{g} && g^d\downarrow\hphantom{g^d} && \tilde g\downarrow\hphantom{\tilde g} \\
B & \stackrel{\tau_B}{\longleftarrow} & B^d  & \stackrel{\pi_B}{\longleftarrow} & \tilde B
\end{array}
\eqno{(\ref{isortiv.standard.form.exmp}.1)}
$$
We see 
during the proof of (\ref{isortiv.standard.form.thm}) 
that the  cover $X\leftarrow X^d$ 
corresponding to the discrete part of the monodromy representation
is
uniquely determined by $g:X\to B$.
%(It is the unique smallest cover that trivializes the
%local system  $R^1g_*\q_X$ over a dense open set.)
By contrast, the $X^d\leftarrow \tilde X$ part is not unique.
Its group of deck transformations is
$G_t\subset \aut^0(F)$, hence Abelian.
It is not even clear that there is a natural ``smallest'' choice of
$X^d\leftarrow \tilde X$.

If $F=A$ is an Abelian variety then 
$g^d:X^d\to B^d$ is  a  Seifert bundle where 
an  orbibundle  $g^s:X^s\to B^s$ is called a {\it Seifert bundle}
if $F=A$ is an Abelian variety and $G$ acts on $A$ by translations.
Note that in this case the $A$-action on $\tilde B\times A$
descends to an $A$-action on $X^s$ and $B^s=X^s/A$.
Thus the reduced structure of every fiber is a
smooth Abelian variety isogenuous to $A$.

\end{defn}

\begin{lem} \label{qt.elliptic.td2=0.lem}
Notation as above. Then
\begin{enumerate}
\item $\pi_X$ and $\tau_X$ are \'etale in codimension 1
(that is, quasi-\'etale),
\item $\pi_X$ and $\tau_X$ are \'etale in codimension 2
if  one of the following holds
\begin{enumerate}
\item  $G$ acts freely on $F$ outside a codimension $\geq 2$ subset or
\item $K_F\sim 0$ and $\Delta_{X/B}=0$.
\end{enumerate}
%\item If, in addition, $\dim A=1$ then $\tau_B$
%is  \'etale in codimension 1.
\end{enumerate}
\end{lem}

Proof.  The first claim  is clear since  both $\rho_F, \rho_B$ are faithful.

Since $\rho_F, \rho_B$ are faithful, $\tau_X$
 fails to be \'etale in codimension 2
iff some $1\neq g\in G$ fixes a  divisor $\tilde D_B\subset \tilde B$
and also a  divisor $D_F\subset  F$. This is excluded by (2.a).

Next we check that (2.b) implies (2.a).
At a general point $p\in D_F$ choose local $g$-equivariant coordinates
$x_1,\dots, x_m$ such that  $D_F=(x_1=0)$.
Thus $\rho_F(g)^*$ acts on $x_1$ non-trivially but it fixes
$x_2,\dots, x_m$.
Let $\omega_0$ be a nonzero section of $\omega_F$.
Locally near $p$ we can write
$$
\omega_0= f\cdot dx_1\wedge \cdots  \wedge dx_m,
$$
thus  $\rho_F(g)^*$ acts on $H^0(F,\omega_F)$
with the same eigenvalue as on $x_1$.

Thus, by (\ref{kodaira.formula}.1),  
the image of $\tilde D_X$   gives a positive contribution
to $\Delta_{X/B}$. This contradicts  $\Delta_{X/B}=0$.\qed

%Finally, if $\dim F=1$ then every non-translation $g\in G$ fixes a
%divisor in $F$ so the previous argument shows that
%$g$ can not fix a divisor on $\tilde B$. 

%Note that (3) can fail if $\dim A>1$ and this makes
%the study of higher dimensional Abelian orbibundles
%more complicated.

\medskip

There are some obvious deformations of $X$ obtained by
deforming $\tilde B$ and $F$ in a family
$\{(\tilde B_t, F_t)\}$
such that the representations $\rho_B,\rho_F$
lift to $\rho_{B,t}:G\to \aut(\tilde B_t)$
and   $\rho_{F,t}:G\to \aut(F_t)$.

In general, not every deformation of $X$ arises this way.
For instance, let $\tilde B$ and $F=A$ be elliptic curves and
$X$ the Kummer surface of $\tilde B\times A$.
The obvious deformations of $X$ form a 2-dimensional family obtained by
deforming $\tilde B$ and $A$.
Thus a  general deformation of $X$ is not 
obtained this way and it is not even 
elliptic.
Even worse, a general  elliptic deformation of $X$
is also not Kummer, thus  not every deformation of 
the morphism $(g:X\to B)$ is obtained by the quotient
construction. 

% However, the next theorem shows that in the cases not covered
% by (\ref{main.char.thm.td2}), every deformation of $X$ arises the obvious way.

\begin{thm} \label{isortiv.standard.form.def.thm}
Let  $g:X\to B$ be a Calabi--Yau orbibundle with general fiber $F$. Assume that
 $X$ has log terminal singularities, $H^2(X, \o_X)=0$,  $\kappa(X)\geq 0$,
$K_F\sim 0$ 
and $\Delta_{X/B}=0$.

Then every flat deformation of $X$ arises from a flat
deformation of $\bigl(\tilde B, F, \rho_B, \rho_F\bigr)$.
\end{thm}

Proof.  Let $L_B$ be an ample line bundle on $B$ and set
$L:=g^*L_B$.

Let $h:{\mathbf X} \to (0\in S)$ be a 
deformation of $X_0\cong X$.
In the sequel we will repeatedly replace $S$ by a smaller analytic (or \'etale)
neighborhood of $0$ if necessary.

 Since $H^2(X, \o_X)=0$, $L$ lifts to a line bundle
${\mathbf L} $ on $ {\mathbf X}$ by   \cite[p.236-16]{FGA}.

Since  $K_F\sim 0$ 
and $\Delta_{X/B}=0$,  (\ref{qt.elliptic.td2=0.lem})
implies that 
$\pi:\tilde X\to X$ is \'etale in codimension 2. Thus, by
\cite[Thm.12]{k-flat}, the cover $\pi$ lifts to
a cover $\Pi: \tilde{\mathbf X}\to {\mathbf X}$. 

Finally we show that the product
decomposition
$\tilde X\cong \tilde B\times F$ lifts to a
product
decomposition
$$
\tilde {\mathbf X}\cong \tilde {\mathbf B}\times_S {\mathbf F}
$$
where $ \tilde {\mathbf B}\to S$ is a  flat deformation of $\tilde B$
and ${\mathbf F}\to S$ is a family of 
Calabi--Yau varieties over $S$. 
After a further \'etale cover of $\tilde F\to F$ we may assume that
$\tilde F\cong Z\times A$ where $H^1(Z,\o_Z)=0$ and $A$ is an Abelian variety.
Set $\hat X:= \tilde B\times Z\times A$; then
$\hat X \to\tilde X$ lifts to a deformation 
$ \widehat {\mathbf X}\to \tilde {\mathbf X}\to S$.

First we use  (\ref{ab.acts.deforms}) to show that the product
decomposition
$\hat X\cong \bigl(\tilde B\times Z\bigr)\times A$ lifts to a
product
decomposition
$$
\widehat {\mathbf X}\cong \widehat {\mathbf {BZ}}\times_S {\mathbf A}
$$
where $ \widehat {\mathbf {BZ}}\to S$ is a  flat deformation of 
$\tilde B\times Z$
and ${\mathbf A}\to S$ is a family of 
Abelian varieties over $S$. 
There are two separate issues here: we have to make sure that
automorphisms of $\hat X$ lift to automorphisms of
$\widehat {\mathbf X} $ and we have to ensure that
 the lifted $A$-action does not get mixed-up with the
possible automorphisms of $\tilde B$. 

The deformation of the product $\tilde B\times Z $
is much easier; we discuss it in 
(\ref{deform.conj.say}). \qed

\begin{say} In some sense, elliptic curves give the only 
examples of  Calabi--Yau orbibundles that have a non-obvious deformation.

Assume that $F$ has no finite \'etale cover $\tilde F$
that can be written as a product  $\tilde F\cong F_1\times E$
where $E$ is an elliptic curve.
We claim that if $H^2(X, \o_X)=0$ then  
every flat deformation of $X$ arises from a flat
deformation of $\bigl(\tilde B, F, \rho_B, \rho_F\bigr)$.

To see this consider the relative Albanese
$$
\begin{array}{ccc}
X^0 & \stackrel{alb}{\longrightarrow} & \alb_{B^0}(X^0)\\
\downarrow && \downarrow \\
B^0 & = & B^0.
\end{array}
$$
If an automorphism of an Abelian variety $A$ fixes a divisor
then that divisor is an Abelian subvariety and
$A$ has an elliptic curve factor up-to isogeny. 
Thus (\ref{qt.elliptic.td2=0.lem}) shows that our arguments
above apply to $\alb_{B^0}(X^0)\to B^0$. 
Hence by taking a suitable cover that is \'etale outside a
codimension $\geq 2$ subset, we can trivialize
$\alb_{B^0}(X^0)\to B^0$. Thus we may assume that
$\alb_{B^0}(X^0)\cong B^0\times \alb(F)$

The Albanese map  $F\to \alb(F)$ is a fiber bundle, thus
after taking a finite cover  $\alb'(F)\to \alb(F)$ we  get
$F'\to F$ such that
$F'\cong F_1\times \alb'(F)$. If $H^1(F_1, \o_{F_1})\neq 0$ then
$\dim \alb(F')>\dim \alb(F)$ and we repeat the above argument.

Thus eventually    we get a cover $\hat X\to X$ that is \'etale outside a
codimension $\geq 2$ subset such that $\hat g: \hat X\to B$ is an
orbibundle with fiber $\hat F$ and there is a morphism
$\hat q: \hat X\to B\times \alb(\hat F)$ which is an orbibundle with
fiber $ G$ with $H^1(G, \o_{G})= 0$.

By \cite[Thm.12]{k-flat}, every deformation $h:{\mathbf X}\to (0\in S)$ of
$X$ lifts to a  deformation $\widehat{\mathbf X}$  of $\hat X$.
Since  $H^2(X, \o_X)=0$, we can lift $L$ to a line bundle
$\mathbf L$ on 
$\mathbf X$, hence to  a line bundle
$\widehat{\mathbf L}$ on 
$\widehat{\mathbf X}$.
We see in (\ref{deform.conj.say}) that $\hat g$
lifts to a morphism  $\hat{\mathbf g}: \widehat{\mathbf X}\to {\mathbf{BA}}$
where $ {\mathbf{BA}}$ is a deformation of $B\times \alb(\hat F)$.
We can use $\hat{\mathbf g}$ to descend $\widehat{\mathbf L}$
to a line bundle $\widehat{\mathbf L}_{BA}$
on $ {\mathbf{BA}}$. 
Now we can use  (\ref{ab.acts.deforms}) to see that 
 $ {\mathbf{BA}}\cong   {\mathbf{B}} \times_S {\mathbf{A}}$
where ${\mathbf{B}}$ is a deformation of $B$ and 
${\mathbf{A}}$ is a deformation of $\alb(\hat F)$. \qed

\end{say}

\begin{lem} \label{ab.acts.deforms}
Let $Y\to S$ be a flat, proper morphism whose fibers are normal
and  $L$  a line bundle on $Y$. Let $0\in S$ be a point
such that
\begin{enumerate}
\item $Y_0$ is not birationally ruled,
\item an Abelian variety $A_0\subset \aut^0(Y_0)$ acts faithfully on $Y_0$,
\item $L_0$ is nef, %$\nu(Y_0, L_0)=\dim Y_0-\dim A_0$ and
$L_0$ is numerically trivial on the $A_0$-orbits
but not numerically trivial on general $A'_0$-orbits
for any  $A_0\subsetneq A'_0\subset\aut^0(Y_0)$.
\end{enumerate}
Then, possibly after shrinking $S$, there is an
Abelian scheme $A\to S$ extending $A_0$ such that
$A$ acts faithfully on $Y$.
\end{lem}

Proof. By \cite[p.217]{MR0258827}
(see also \cite[p.392]{Kollar85})  possibly after shrinking $S$,
$g^a:\aut^0(Y/S)\to S$ is a smooth Abelian scheme,
where $\aut^0(Y/S) $ denotes the identity component of the
automorphism scheme  $\aut(Y/S) $. 
Working \'etale locally, we may assume that there
is a section $Z\subset Y$. Acting on $Z$ gives a morphism
$\rho_Z:\aut^0(Y/S)\to Y$. Then $\rho_Z^*L$ is a nef line bundle
on $\aut^0(Y/S)$.  
The kernel of the cup-product map
$$
c_1\bigl(\rho_Z^*L\bigr) :R^1g^a_*\ \q\to R^3g^a_*\ \q
$$
is a variation of sub-Hodge structures, hence it
corresponds to a smooth Abelian subfamily
$A\subset \aut^0(Y/S)$. By (3), this is the required
extension of $A_0$.

The quotient then exists by \cite{MR0164973}. 
\qed
\medskip

We will also need to understand the class group of an orbibundle.

\begin{say}[Divisors on orbibundles]
\label{isortiv.divisors.say}
We use the notation of (\ref{isortiv.standard.form.exmp})
and of (\ref{albanese.defns}).

 By \cite[5.3]{2011arXiv1104.1861B}, 
(see also \cite{hul-klo, 2010arXiv1008.2018C} for the elliptic case)
the  class group of $\tilde B\times F $ is
$$
\cl\bigl(\tilde B\times F\bigr)=\cl(\tilde B)+\cl(F)+
\Hom\bigl(\albr(\tilde B), \pic^0(F)\bigr).
\eqno{(\ref{isortiv.divisors.say}.1)}
$$
This comes with a natural $G$-action and, up-to torsion,
the class group of the quotient is
$$
\cl(B)+\cl(F)^G+\Hom\bigl(\albr(\tilde B), \pic^0(F)\bigr)^G.
\eqno{(\ref{isortiv.divisors.say}.2)}
$$
If $\tilde B$ has rational singularities then
$\albr(\tilde  B)=\alb(\tilde B)$ and then the extra component 
$\Hom\bigl(\alb(\tilde  B), \pic^0(F)\bigr)$ gives
$\q$-Cartier divisors.

We will use the following variant of this.

\medskip

{\it Claim} \ref{isortiv.divisors.say}.3.
Let $g:X\to B$ be an orbibundle such that $X$ has  log terminal  singularities.
Then the natural map
$$
\cl(B)/\pic(B)+\bigl(\cl(F)/\pic(F)\bigr)^G\to \cl(X)/\pic(X)
$$
is an isomorphism modulo torsion.
In particular, 
if  $B$ and $F$ are  $\q$-factorial then so is $X$. 
\medskip

%{\it Claim} \ref{isortiv.divisors.say}.4.
%Let $g:X\to B$ be an orbibundle.
%Assume that $\Delta_{X/B}=0$ and
%$B$  and $F$ are  $\q$-factorial, log terminal.
%Then $X$ also has $\q$-factorial singularities.

Proof. By (\ref{qt.elliptic.td2=0.lem}),
 $\tau_X:X^d\to X$ is \'etale in codimension 1,
hence $X^d$ also has log terminal singularities.
As noted in (\ref{Calabi--Yau.fib.defn}), this implies that
$B^d$ has rational singularities.

Let us now study more carefully the right hand side of
(\ref{isortiv.divisors.say}.2). Let $G_t\subset G$
denote the subgroup of translations. Then
$$
\Hom\bigl(\albr(\tilde B), \pic^0(F)\bigr)^G\subset
\Hom\bigl(\albr(\tilde B), \pic^0(F)\bigr)^{G_t}.
$$
Since translations act trivially on $\pic^0(F) $,
the latter can be identified (up-to torsion) as
$$
\begin{array}{rcl}
\Hom\bigl(\albr(\tilde B), \pic^0(F)\bigr)^{G_t}\otimes\q &\cong&
\Hom\bigl(\albr(\tilde B)^{G_t}, \pic^0(F)\bigr)\otimes\q\\
&\cong&
\Hom\bigl(\albr( B^d), \pic^0(F)\bigr)\otimes\q\\
&\cong&\Hom\bigl(\alb( B^d), \pic^0(F)\bigr)\otimes\q.
\end{array}
$$
Thus this extra term gives only $\q$-Cartier divisors on $X^d$
and hence also on $X$. \qed

\medskip
The following local example shows that it
 is not enough to assume that $B$ has rational singularities.
Set $\tilde B=(u^3+v^3+w^3=0)\subset \a^3$ and
$E=(x^3+y^3+z^3=0)\subset \p^2$. On both factors, $\z/3$ acts
by weights $(0,0,1)$. Then $B=\tilde B/\frac13(0,0,1)\cong \a^2$
is even smooth but
$$
X=\tilde B\times E/\tfrac13(0,0,1)\times (0,0,1)
$$
is not $\q$-factorial. For instance, the closure of the graph of the
natural projection $\tilde B\map E$
gives a non-$\q$-Cartier divisor on $X$.
\end{say}

\begin{defn}[Albanese varieties]\label{albanese.defns}
For a smooth  projective variety $V$ let 
 $\alb(V)$ denote the  Albanese variety, that is,
the target of the universal morphism from $V$ to
an Abelian variety.  (See  \cite[Sec.I.13]{bpv} 
or \cite[p.236-16]{FGA} for  introductions.)

There are 2 ways to generalize this concept to
normal varieties. 

The above definition yields  what we again call the
 {\it Albanese variety} $\alb(V)$.
Alternatively,  
 the {\it rational Albanese} variety $\albr(V)$
is defined as the target of the  universal rational map from $V$ to
an Abelian variety.
One can identify $\albr(V)=\alb(V')$ where $V'\to V$ is any
resolution of singularities. 

It is easy to see that if $V$ has log terminal (more generally rational)
 singularities
then  $\albr(V)=\alb(V)$.
\end{defn}

\section{Generically isotrivial  Calabi--Yau fiber spaces}

In this section we prove that all
generically isotrivial  Calabi--Yau fiber spaces
are  essentially   Calabi--Yau orbibundles.

\begin{thm} \label{isortiv.standard.form.thm}
Let $g:X\to B$ be a  projective, generically isotrivial,
Calabi--Yau fiber space. 
\begin{enumerate}
\item There is a unique  Calabi--Yau orbibundle
$\bigl(g^{\rm orb}:X^{\rm orb}\to B\bigr)$
that is  birational to $g:X\to B$.
\item $X$ is isomorphic to  $X^{\rm orb}$ if  the following hold
\begin{enumerate}
\item $X$ is   $\q$-factorial and log terminal,
\item  $g:X\to B$ is relatively minimal and has no exceptional divisors,
\item $B$ is  $\q$-factorial.
\end{enumerate}
\end{enumerate}
\end{thm}

Proof. Let $B^0\subset B$ be a Zariski open subset over which
$X^0\to B^0$ is isotrivial with general fiber $F$. 
This gives a well-defined representation
$$
\rho:\pi_1(B^0)\to \aut(F)/\aut^0(F).
$$
Let $B^{(d,0)}\to B^0$ be the corresponding \'etale, Galois cover
with group $G_d$ 
and $ B^d\to B$ its extension to a (usually ramified)
Galois cover of $B$ with group $G_d$.
This gives the well-defined  cover in
(\ref{isortiv.standard.form.exmp}.1).

The trivialization of the translation part is more subtle
and it depends on additional choices.

A general $ \aut^0(F)$-orbit $A_F\subset F$ defines
an isotrivial Abelian family  $X^{(d,0)}\supset A^{(d,0)}_X\to B^{(d,0)}$.
By assumption there is a $g$-ample line bundle $L$ on $X$.
It pulls back to a relatively ample line bundle $L_A$ on $A^{(d,0)}_X$.
We may assume that its degree on the general fiber is at least 3.
Let $T^{(d,0)}\subset A_X^{(d,0)}$ be the 
subscheme as in (\ref{multsect.constr.say}). 
Since $L_A$ is $G_d$-invariant, $T^{(d,0)}$ is $G_d$-equivariant
hence
it defines a monodromy representation of
$\pi_1(B^0)\to \aut(F)$; let $G$ denote its image.

Let $\tilde B^0\to B^0$ be the corresponding \'etale, Galois cover
with group $G$ 
and $\tilde B\to B$ its extension to a (usually ramified)
Galois cover of $B$ with group $G$. 

By pull-back we obtain an isotrivial,
Abelian fiber space $\tilde A_X^0\to \tilde B^0$
with a trivialization of the $m$-torsion points. For $m\geq 3$
this implies that $\tilde A_X^0\cong \tilde B^0\times A$. 
(This is quite elementary, cf.\ \cite[p.513]{acg2}.)
Thus the same pull-back also trivializes $X^0\to B^0$.
 We can compactify  $\tilde  X^0$ as
$\tilde X:= \tilde B\times A$. 

The $G$-action on $\tilde A_X^0\cong \tilde B^0\times A$
can be given as
$$
g: (\tilde b, a)\mapsto 
\bigl(\rho_B(g)\cdot \tilde b,\rho_{A,\tilde b}(g)\cdot a\bigr).
$$
Note that $\rho_{A,\tilde b} $ preserves the $m$-torsion points
and the automorphisms of an Abelian torsor  that preserve
any finite nonempty set  form a discrete group.
Thus in fact $\rho_{A,\tilde b} $ is independent of $\tilde b$
and hence the $G$-action on $\tilde X$ is given by
$$
g: (\tilde b, a)\mapsto 
\bigl(\rho_B(g)\cdot \tilde b,\rho_{A}(g)\cdot a\bigr)
$$
for some isomorphism $\rho_B:G\cong \gal\bigl(\tilde B/B\bigr)$
and homomorphism $\rho_B:G\to \aut(F)$. 
 We can replace $\tilde B $
by $\tilde B/\ker \rho_B$, hence we may assume that 
 $\rho_B$ is faithful.
By construction $X$ is birational to 
$X^{\rm orb} :=\tilde X/G$. 

In general, birational maps between relatively minimal models
are very special. First there are divisorial contractions along which
the canonical class is trivial. In our case these are excluded by
(2.a). In the non-$\q$-factorial case there could be small  contractions,
but   $X^{\rm orb}$ is also  $\q$-factorial  by 
(\ref{isortiv.divisors.say}.3).

Finally there can be flops, but  the orbibundle does not have
any suitable extremal rays by  (\ref{rho=1.unique.model.lem}).
Thus $X$ is isomorphic to $X^{\rm orb}$
if the conditions (2.a--c) hold.
 \qed

\begin{say}[Multisections of Abelian families]\label{multsect.constr.say}
 Let $E$ be a smooth projective curve of genus 1
and $L$ a line bundle of degree $m$ on $E$. 
If $m=1$ then $L$ has a unique section, thus we can associate a
point $p\in E$ to $L$.  If $m\geq 2$, then sections define a
linear equivalence class $|L|$ of $m$ points. If we fix a point
$0\in E$ to be the origin, then we can add these $m$ points together and 
get a well defined point of $E$ associated to $L$. This, however,
depends on the choice of the origin.

To get something invariant, lets us look at the points
$p\in E$ such that $m\cdot p\in |L|$. There are $m^2$ such points,
together forming a translate of the subgroup of $m$-torsion
points. This construction also works in families.

Let $g:X\to B$ be a smooth, projective morphism
whose fibers $E_b$ are   curves of genus 1.
Let $L$ be a line bundle on $X$ that has degree $m$ on each fiber.
 Then there is a closed subscheme
$T\subset X$ such that $g|_T:T\to B$ is \'etale of degree $m^2$ and
every fiber $T_b\subset E_b$ is a 
 a translate
of the subgroup of $m$-torsion points.

There is a similar construction for higher dimensional Abelian varieties.
For clarity, I say {\it Abelian torsor} when talking about an Abelian variety 
without a specified origin.

Thus let $A$ be an Abelian torsor  of dimension $d$
and $L$  an ample line bundle on $A$. 
It has a first Chern class $\tilde c_1(L)$ in the Chow group
and we get $\tilde c_1(L)^d$ as an element of the Chow group of 0-cycles.
(It is important to use the Chow group, the Chern class in cohomology is not
sufficient.)
Let its degree be $m$. 

Fix a base point $0\in A$. This defines a map
from the Chow group of 0-cycles to $(A,0)$; let
$\alpha\bigl(\tilde c_1(L)^d\bigr)$ denote the image.

Finally let $T\subset A$ be the set of points $t\in A$ such that
$m\cdot t=\alpha\bigl(\tilde c_1(L)^d\bigr)$. This $T$ is  a translate
of the subgroup of $m$-torsion points.
As before, the key point is that $T$ is independent of the choice of the 
base point $0\in A$. Indeed, if we change $0$ be a translation by
$c\in A$ then $\alpha\bigl(\tilde c_1(L)^d\bigr)$ is changed by
translation by
$m\cdot c$ so $T$ is changed by
translation by $c$. 

Furthermore, if $(A_b, L_b)$ is a family of polarized  Abelian torsors
that varies analytically (or algebraically) with $b$ then
$T_b\subset A_b$ is a family of subschemes that also vary
analytically (or algebraically) with $b$. 
Thus we obtain the following.

Let $g:X\to B$ be a  smooth,  projective morphism whose
fibers are 
Abelian torsors. Then there is a closed subscheme
$T\subset X$ such that $g|_T:T\to B$ is \'etale and
every fiber $T_b\subset A_b$ is a 
 a translate
of the subgroup of $m$-torsion points 
(where $\deg T/B=m^{2d}$).
\end{say}

\begin{lem} \label{rho=1.unique.model.lem}
Let $g_i:X_i\to B$ be  projective fiber spaces
and $\phi:X_1\map X_2$ a rational map.
Assume the following.
\begin{enumerate}
\item There are no $g_i$-exceptional divisors.
\item A divisor on $X_2$ is $\q$-Cartier iff its restriction 
to the generic fiber of $g_2$ is $\q$-Cartier.
(This holds trivially if $X_2$ is $\q$-factorial.)
\item Every curve
$C\subset X_2$ contracted by $g_2$ is $\q$-homologous to a
curve in a general fiber. 
\item $\phi$ induces an isomorphism of the generic fibers of the $g_i$.
\item There  are closed   subsets $Z_i\subset X_i$
 such that
$\codim_{X_i}Z_i\geq 2$
and $\phi$ induces an isomorphism  $X_1\setminus Z_1\cong X_2\setminus Z_2$.
\end{enumerate}
Then $\phi$ is an isomorphism. 
\end{lem}

Proof. Let $H_1\subset X_1$ be a $g_1$-ample divisor and
$H_2\subset X_2$ its birational transform.
It follows from assumption (2) and (4) that $H_2$ is $\q$-Cartier and 
from (3) that it is $g_2$-ample.
Thus  (5) and a lemma of Matsusaka--Mumford \cite{ma-mu}
implies that $\phi$ is an isomorphism.
(See \cite[5.6]{ksc} or \cite[Exrc.75]{k-exrc}  for the variant used here.)
  \qed

\begin{say}[F-theory examples] \label{isotriv.section.few}

Let $X$ be a smooth, projective variety and $g:X\to B$  a 
relatively minimal  elliptic fiber space
with a section $\sigma:B\to X$.
Since $X$ is smooth, so is $B$.

Assume that $\Delta_{X/S}=0$.
Then, by (\ref{qt.elliptic.td2=0.lem}),
it can have only multiple smooth fibers at
codimension 1 points, but then  the section shows that
there are no multiple fibers.
Thus there is an an open subset $B^0\subset B$ such that
$\codim_B(B\setminus B^0)\geq 2$ and 
$X^0\to X$ is a fiber bundle with fiber a pointed 
elliptic curve $(E,0)$. Thus $X^0$ is given by the data
$$
\bigl(B^0, E, \rho:\pi_1(B^0)\to \aut(E,0)\bigr).
$$
Note that $\pi_1(B^0)=\pi_1(B)$
since $B$ is smooth and $\codim_B(B\setminus B^0)\geq 2$. 
Thus $X$ is birational to a fiber bundle $g':X'\to B$
given by the data
$$
\bigl(B, E, \rho:\pi_1(B)\to \aut(E,0)\bigr).
$$
All the fibers of $g'$ are elliptic curves but
the exceptional locus of a flip or a flop is always
covered by rational curves (cf.\ \cite[VI.1.10]{rc-book}).
Thus in fact $X\cong X'$ hence $g:X\to B$ is a locally trivial fiber bundle.
The image of the monodromy representation
$\rho:\pi_1(B)\to \aut(E,0) $ is usually $\z/2$, but 
for elliptic curves with extra automorphisms 
it can also be $\z/3, \z/4$ or $\z/6$.

It is easy to write down  examples where $K_X\sim 0$ and
$H^i(X,\o_X)=0$ for $0<i<\dim X$. However,  $\pi_1(X)$ is always infinite,
so such an $X$ can not be a ``true'' Calabi--Yau manifold.

By (\ref{isortiv.standard.form.def.thm}),
 if $H^2(X, \o_X)=0$ then every small deformation of
$X$ is obtained by deforming $B$ and, if the image of $\rho$ is $\z/2$,
also deforming $E$.
\end{say}

\section{Examples}

The first example is an elliptic Calabi--Yau surface  with
quotient singularities 
that has a flat smoothing
which is neither  Calabi--Yau nor  elliptic.

\begin{exmp} \label{non.CY.def.exmp}
We start with  a surface 
 $S^*_F$  which is
the  quotient of the square of the Fermat cubic curve by $\z/3$:
$$
S^*_F\cong \bigl(u_1^3=v_1^3+w_1^3\bigr)\times \bigl(u_2^3=v_2^3+w_2^3\bigr)/
\tfrac13(1,0,0;1,0,0).
$$
To describe the deformation, we need a different representation of it.

In $\p^3$ consider two lines
$L_1=(x_0=x_1=0)$ and $L_2=(x_2=x_3=0)$. 
The linear system $\bigl|\o_{\p^2}(2)(-L_1-L_2)\bigr|$
is spanned by the 4 reducible quadrics  $x_ix_j$ for $i\in \{0,1\}$ and
$j\in \{2,3\}$. They satisfy a relation  
$(x_0x_2)(x_1x_3)=(x_0x_3)(x_1x_2)$. Thus we get a morphism
$$
\pi: B_{L_1+L_2}\p^3\to \p^1\times \p^1
$$ 
which is a $\p^1$-bundle whose fibers are the birational transforms
of lines that intersect both of the $L_i$.

Let $S\subset \p^3$ be a cubic surface such that
${\mathbf p}:=S\cap (L_1+L_2)$ is 6 distinct points.
Then we get
$\pi_S: B_{\mathbf p}S\to  \p^1\times \p^1$. 

In general, none of the lines connecting 2 points of ${\mathbf p} $
is contained in $S$. Thus in this case $\pi_S$ is a finite triple cover.

Both of the lines  $L_i$ determine an elliptic pencil on 
$B_{\mathbf p}S $ but if we move the 6 points 
${\mathbf p}$ into general position, we lose both
elliptic pencils.

At the other extreme we have the Fermat-type surface
$$
S_F:=\bigl(x_0^3+x_1^3=x_2^3+x_3^3\bigr)\subset \p^3.
$$
We can factor both sides and write its equation as
$m_1m_2m_3=n_1n_2n_3$.
The 9 lines  $L_{ij}:=(m_i=n_j=0)$ are all contained in
$S_F$. Let $L'_{ij}\subset  B_{\mathbf p}S_F$ denote their birational
transforms. Then the self-intersections
$\bigl(L'_{ij}\cdot L'_{ij}\bigr)$ equal $-3$ and $\pi_{S_F}$
contracts these 9 curves $L'_{ij}$. 
Thus the Stein factorization of $\pi_{S_F}$
gives a triple cover  $S^*_F\to \p^1\times \p^1$
and $S^*_F$ has 9 singular points of type
$\a^2/\frac13(1,1)$.
We see furthermore that
$$
-3K_{S_F}\sim \tsum_{ij} L_{ij}\qtq{and}
-3K_{B_{\mathbf P}S_F}\sim \tsum_{ij} L'_{ij}.
$$
Thus $-3 K_{S^*_F}\sim 0$.

To see that this is the same  $S^*_F$, note that 
the morphism of  the original $S^*_F$ to $\p^1\times \p^1$ is given by
$$
(u_1{:}v_1{:}w_1)\times (u_2{:}v_2{:}w_2)\mapsto
(v_1{:}w_1)\times (v_2{:}w_2)
$$
and the rational map to the cubic surface is given by
$$
(u_1{:}v_1{:}w_1)\times (u_2{:}v_2{:}w_2)\mapsto
\bigl(v_2u_1u_2^2{:} u_1u_2^2{:} v_1u_2^3{:} u_2^3\bigr).
$$

Varying $S$ gives a flat deformation whose central fiber is
$S^*_F$, a surface with quotient singularities and
torsion canonical class and whose general fiber is
a cubic surface blown up at 6 general points, hence rational
and without elliptic pencils.
\end{exmp}

The next example gives local models of
generically isotrivial elliptic orbibundles
that have a crepant resolution. 

\begin{exmp} \label{small.res.emp}
Let $Z\subset \p^N$ be an anticanonically embedded
Fano variety and $X\subset \a^{N+1}_{\mathbf x}$ the cone over $Z$. 
Let $0\in E$ be an elliptic curve with a marked point. 
Consider the elliptic fiber space
$$
Y:=X\times E/(-1,-1) \to X/(-1).
$$
We claim that $Y$ has a crepant resolution.

First we blow up the vertex of $X$. We get
$B_0X\to X$ with exceptional divisor $F\cong Z$.
Note further that $B_0X\to X$ is crepant.
The involution lifts to 
$B_0X\times E/(-1,-1)$. The fixed point set of this action is
$F\times \{0\}$; a smooth subvariety of codimension 2.
Thus $B_0X\times E/(-1,-1)$ is resolved by
blowing up the singular locus.
\end{exmp}

The next example shows that for surfaces with normal crossing singularities,
a deformation may lose the elliptic structure.

\begin{exmp}\label{rtle.elliptic.bad.exmp}
 Let $S\subset \p^1\times \p^2$
be a smooth surface of bi-degree $(1,3)$. The first projection
$\pi:S\to \p^1$ is an elliptic fiber space.
The other projection $\tau:S\to \p^2$
exhibits it as the blow-up of $\p^2$ at 9 base points
of an elliptic pencil.
Let $F_1,\dots, F_9\subset S$ denote the 9 exceptional curves.
Thus $S$ is an elliptic $dP_9$.
In particular,  specifying $\pi:S\to \p^1$ plus a fiber of $\pi$
is equivalent to a pair 
 $\bigl(E\subset \p^2\bigr)$ 
plus 9 points  $P_1,\dots, P_9\in E$ 
such that $P_1+\cdots + P_9\sim \o_{\p^2}(3)|_E$.
The elliptic pencils are given by 
$\pi^*\o_{\p^1}(1)\cong \tau^*\o_{\p^2}(3)(-F_1-\cdots-F_9)$.

Let us now vary the points on $E$
in a family $P_i(t): t\in \c$.
The line bundle giving the
elliptic pencil deforms as
$\tau^*\o_{\p^2}(3)(-F_1(t)-\cdots-F_9(t))$
but the elliptic pencil deforms only if
$P_1(t)+\cdots + P_9(t)\sim \o_{\p^2}(3)|_E$ holds for every $t$.

Let $X\subset \p^2\times \p^2$
be a smooth 3-fold of bi-degree $(1,3)$. The first projection
$\pi:X\to \p^2$ is an elliptic fiber space.

If $C\subset \p^2$ is a conic, its preimage
$X_C\to C$ is an elliptic K3 surface.
If $C$ is general then $X_C$ is smooth.

If $C=L_1\cup L_2$ is a pair of  general lines then
$X_C=S_1\cup S_2$ is a singular K3 surface which is a
union of 2 smooth $dP_9$ that intersect along a smooth elliptic curve $E$.

We can thus think of $X_C$ as obtained from
two pairs $\bigl(E^i\subset \p^2\bigr)$ ($i=1,2$) 
with an isomorphism $\phi:E^1\to E^2$
by blowing up 9 points $P^i_j\subset E^i$ ($j=1,\dots, 9$)
 and gluing the resulting surfaces
along the birational transforms of $E^1$ and $E^2$.

Let us now vary the  points on both curves
$P^1_i(t)$ and $P^2_i(t)$. 
We get two families $S_1(t), S_2(t)$ 
and this induces a deformation
$X_C(t)=S_1(t)\cup S_2(t)$.

Although the line bundle $\pi^*\o_C(1)$ giving the
elliptic pencil $X_C\to C$ deforms
on both of the $S_i(t)$, in general we do not get
a line bundle on $X_C(t)$ unless
$$
P^1_1(t)+\cdots + P^1_9(t)\sim 
\phi^*\bigl(P^2_1(t)+\cdots + P^2_9(t)\bigr)
$$
holds for every $t$.
We can thus arrange that 
$\pi^*\o_C(1)$ deforms along $X_C(t)$ but
we lose the  elliptic pencil.
\end{exmp}

\section{General conjectures}

A straightforward generalization of Conjecture \ref{ques4}
is the following, cf.\  \cite{oguiso93} and \cite[Lect.10]{MP97}.

\begin{conj}[Strong abundance for Calabi--Yau manifolds] 
\label{ques3}
Let $X$ be a Calabi--Yau manifold and  $L\in H^2(X, \q)$ a $(1,1)$-class
such that $(L\cdot C)\geq 0$ for every algebraic curve $C\subset X$.
Then there is a unique morphism with connected fibers
$g:X\to B$ onto a normal variety $B$ and  an ample 
$L_B\in H^2(B, \q)$  such that
$L= g^*L_B$. 
\end{conj}

The usual abundance conjecture assumes that $L$ is effective,
but this may not be necessary.

One expects that %Conjecture 
(\ref{ques3}) gets harder as the
 dimension of $B$ decreases. The easiest case, when
$\dim B=\dim X-1$ corresponds to  Questions \ref{ques1}--\ref{ques2}.

From the point of view of higher dimensional birational geometry,
it is natural to consider a more general setting.

%morphisms $g:X\to B$ such that $K_X$ is trivial on the fibers of $g$. 
%In particular, the general fiber is a  Calabi--Yau variety.
%At least conjecturally, everything seems to generalize to
%this setting, at least when $X$ has nonnegative Kodaira dimension.

A {\it log Calabi--Yau fiber space} is a
proper morphisms
with connected fibers  $g:(X,\Delta)\to B$
onto a normal variety
where $(X,\Delta) $ is klt (or possibly lc) and
 $(K_X+\Delta)|_{X_g}\simq 0$ where $X_g\subset X$ is a general fiber.

%The aim of this section is to  reformulate
%Conjecture \ref{ques3}  for
%log Calabi--Yau fiber spaces and to relate it to some of the
%standard conjectures of the Minimal Model Program.

%The most important variant is (\ref{charact.conj}),
%which gives an abstract characterization of  log Calabi--Yau fiber spaces.
%It  would imply that a small deformation of a
%log Calabi--Yau fiber space is again a log Calabi--Yau fiber space (\ref{deform.conj}).

Let $(X,\Delta) $ be a proper klt pair such that 
$K_X+\Delta$ is nef and
$g:(X,\Delta)\to B$  a relatively minimal Calabi--Yau fiber space.
Let $L_B$ be an ample $\q$-divisor on $B$
and set $L:=g^*L_B$. Then 
$L-\epsilon(K_X+\Delta)$ is nef for $0\leq \epsilon\ll 1$.
The converse fails in some rather simple cases, for instance
when $X=B\times E$ for an elliptic curve $E$ and 
we twist $L$ by a degree zero non-torsion line bundle on $E$.

It is natural to expect that the above are essentially the only exceptions.

\begin{conj} \label{charact.conj} 
Let $(X,\Delta) $ be a proper klt pair such that 
$K_X+\Delta$ is nef and $H^1(X,\o_{X})=0$.
Let $L$ be a Cartier divisor on $X$
such that   $L-\epsilon(K_X+\Delta)$ is nef for $0\leq \epsilon\ll 1$.

Then there is a relatively minimal log Calabi--Yau fiber space
structure $g:(X,\Delta)\to B$ and an ample $\q$-divisor
$L_B$  on $B$ such that $L\simq g^*L_B$.
\end{conj}

If  $L-\epsilon(K_X+\Delta)$ is effective then (\ref{charact.conj})
is implied by the Abundance Conjecture.
Note also that (\ref{rtle.elliptic.bad.exmp}) shows that (\ref{charact.conj})
fails if $(X,\Delta) $ is log canonical.

% As in the proof of  (\ref{main.char.thm.cor.q1}), we see that
% (\ref{charact.conj}) implies the following.
%Proof of (\ref{charact.conj}) $\Rightarrow$ (\ref{deform.conj}).

\begin{conj} \label{deform.conj} 
Let  $g_0:(X_0,\Delta_0)\to B_0$ be
 a log Calabi--Yau fiber space where   $(X_0,\Delta_0) $ is 
 a proper klt pair and
$H^2(X_0,\o_{X_0})=0$. 

Let $(X,\Delta) $ be a  klt pair and
$h:(X, \Delta)\to (0\in S)$  a flat proper morphism
whose central fiber is $(X_0,\Delta_0) $.

Then, after passing to an analytic or \'etale
neighborhood of $0\in S$,  there is a proper, flat morphism $B\to (0\in S)$
whose central fiber is $B_0$ such that $g_0$ extends to a
log Calabi--Yau fiber space  $g:(X,\Delta)\to B$.
\end{conj}

\begin{say} \label{deform.conj.say}
Although (\ref{deform.conj}) looks much more general
than (\ref{thm.ques1.general}), it seems that 
Abelian fibrations comprise the only unknown case.

Indeed, let $X_0, B_0$ be projective varieties with
rational singularities and $g_0:X_0\to B_0$  a  morphism 
with connected general fiber $F_0$.
Assume that $H^1\bigl(F_0, \o_{F_0}\bigr)=0$. 
Then $R^1(g_0)_* \o_{X_0}$ is a torsion sheaf. On the other hand,
it is reflexive by \cite[7.8]{k-hdi1}. Thus 
$R^1(g_0)_* \o_{X_0}=0$.

Let $L_{B_0}$ be a sufficiently ample line bundle on $B_0$
and set $L_0:=g_0^*L_{B_0}$. Then 
$H^1\bigl(X_0, L_0\bigr)=0$ by (\ref{22.say}.1). Thus, if 
$h:X\to (0\in S)$  is a deformation of $X_0$ such that
$L_0$ lifts to a line bundle $L$ on $X$
then every section of $L_0$ lifts to a section of $L$
(after passing to an analytic or \'etale
neighborhood of $0\in S$).
Thus (\ref{deform.conj}) holds in this case.

Furthermore, the method of (\ref{smooth.horik.thm}) suggests that
the most difficult case is Abelian pencils over $\p^1$.

Note also that it is easy to write down examples
of Abelian Calabi--Yau fiber spaces $f:X\to B=\p^1$
such that $\Hom_B\bigl(\Omega_B,R^1f_*\o_X\bigr)\neq 0$, thus
(\ref{horik.thm}) does not seem to be sufficient to prove
(\ref{deform.conj}). 
\end{say}

\begin{say}[Log elliptic fiber spaces]\label{logell.conj.say}
As before,  $g:(X,\Delta)\to B$ is a log elliptic
fiber space iff $\bigl(L^{\dim X}\bigr)=0$ but 
$\bigl(L^{\dim X-1}\bigr)\neq 0$. There are 3 cases to consider.
\begin{enumerate}
\item If $\bigl(L^{\dim X-1}\cdot \Delta\bigr)> 0$ then
Riemann--Roch shows that $h^0(X, L^m)$ grows like
$m^{\dim X-1}$ and we get (\ref{charact.conj})
as in (\ref{main.char.thm.td2}).
In this case the general fiber of $g$ is  $F\cong \p^1$
and $\bigl(F\cdot \Delta\bigr)=2$. 
\item If $\bigl(L^{\dim X-1}\cdot \Delta\bigr)= 0$
but $\bigl(L^{n-2}\cdot \operatorname{td}_2(X)\bigr)>0$ then
the proof of (\ref{main.char.thm.td2}) works with minor changes.
\item The hard and unresolved case is again when
$\bigl(L^{\dim X-1}\cdot \Delta\bigr)= 0$
and $\bigl(L^{n-2}\cdot \operatorname{td}_2(X)\bigr)=0$, 
 so $\chi(X, L^m)=O\bigl(m^{\dim X-3}\bigr)$. 
\end{enumerate}
% In the last case $g:(X,\Delta)\to B$ is birational to a
%  Calabi--Yau orbibundle
% $g^{\rm orb}:\bigl(X^{\rm orb},\Delta^{\rm orb}\bigr)\to B$
%  by (\ref{isortiv.standard.form.thm}) and every deformation of 
% such an orbibundle is again log elliptic by 
% (\ref{isortiv.standard.form.def.thm}). 

% Note that  we have not fully proved that every deformation of 
%  a log elliptic
% fiber space $g:(X,\Delta)\to B$ is again log elliptic. 
% If ${\mathbf X}\to S$ is a family of smooth 3-folds then
% the MMP for the family is a also an MMP on each fiber
% by \cite{km-flips}. Thus, in this case, the deformation
% theory of a minimal model encompasses the  deformation
% theory of every smooth model. This  fails for
% smooth 4-folds or for singular 3-folds. Thus I do not have a
% complete answer
% to Question \ref{ques1}.  However, as we noted in
% (\ref{various.say}),  our results are sufficient to yield a positive answer
% in the cases that arise in F-theory, that is
% when the elliptic fibration has a section. 
\end{say}

 \begin{ack}
I thank G.~Di~Cerbo, 
R.~Donagi, O.~Fujino, A.~Langer, R.~Lazarsfeld, K.~Oguiso, Y.-C.~Tu and C.~Xu 
for helpful discussions, comments and references.
Partial financial support   was provided  by  the NSF under grant number 
DMS-07-58275.
\end{ack}

%\bibliography{refs-main/refs}
\def\cprime{$'$} \def\cprime{$'$} \def\cprime{$'$} \def\cprime{$'$}
  \def\cprime{$'$} \def\cprime{$'$} \def\dbar{\leavevmode\hbox to
  0pt{\hskip.2ex \accent"16\hss}d} \def\cprime{$'$} \def\cprime{$'$}
  \def\polhk#1{\setbox0=\hbox{#1}{\ooalign{\hidewidth
  \lower1.5ex\hbox{`}\hidewidth\crcr\unhbox0}}} \def\cprime{$'$}
  \def\cprime{$'$} \def\cprime{$'$} \def\cprime{$'$}
  \def\polhk#1{\setbox0=\hbox{#1}{\ooalign{\hidewidth
  \lower1.5ex\hbox{`}\hidewidth\crcr\unhbox0}}} \def\cdprime{$''$}
  \def\cprime{$'$} \def\cprime{$'$} \def\cprime{$'$} \def\cprime{$'$}
\providecommand{\bysame}{\leavevmode\hbox to3em{\hrulefill}\thinspace}
\providecommand{\MR}{\relax\ifhmode\unskip\space\fi MR }
% \MRhref is called by the amsart/book/proc definition of \MR.
\providecommand{\MRhref}[2]{%
  \href{http://www.ams.org/mathscinet-getitem?mr=#1}{#2}
}
\providecommand{\href}[2]{#2}

\vskip1cm

\noindent Princeton University, Princeton NJ 08544-1000

{\begin{verbatim}kollar@math.princeton.edu\end{verbatim}}

\end{document}